\documentclass[11pt]{amsart}
\usepackage[T1]{fontenc}
\usepackage[utf8]{inputenc}
\usepackage[english]{babel}
\usepackage{lmodern,amsmath,amssymb,amsthm,mathtools,microtype}
\usepackage{enumitem,array,booktabs}
\usepackage{tikz}
\usetikzlibrary{arrows.meta}
\usepackage[a4paper,margin=2.8cm]{geometry}
\usepackage[hypertexnames=false,hidelinks]{hyperref}
\numberwithin{equation}{section}
\makeatletter
\let\shapironormalsize\normalsize
\renewcommand{\normalsize}{\shapironormalsize
  \setlength{\abovedisplayskip}{10pt plus 2pt minus 2pt}\setlength{\belowdisplayskip}{10pt plus 2pt minus 2pt}\setlength{\abovedisplayshortskip}{7pt plus 2pt minus 1pt}\setlength{\belowdisplayshortskip}{8pt plus 2pt minus 1pt}\setlength{\jot}{4pt}}
\renewcommand{\section}{\@startsection{section}{1}{\z@}{-20pt plus 4pt minus 2pt}{10pt plus 2pt minus 1pt}{\normalfont\large\bfseries\raggedright}}
\renewcommand{\subsection}{\@startsection{subsection}{2}{\z@}{-15pt plus 3pt minus 2pt}{7pt plus 1pt minus 1pt}{\normalfont\normalsize\bfseries\raggedright}}
\renewcommand{\subsubsection}{\@startsection{subsubsection}{3}{\z@}{-12pt plus 2pt minus 1pt}{6pt plus 1pt minus 1pt}{\normalfont\normalsize\bfseries\raggedright}}
\renewcommand{\paragraph}{\@startsection{paragraph}{4}{\z@}{-10pt plus 2pt minus 1pt}{5pt plus 1pt minus 1pt}{\normalfont\normalsize\itshape\raggedright}}
\makeatother
\AtBeginDocument{\normalsize}
\theoremstyle{plain}
\newtheorem{theorem}{Theorem}[section]
\newtheorem{lemma}[theorem]{Lemma}
\newtheorem{proposition}[theorem]{Proposition}
\newtheorem{corollary}[theorem]{Corollary}

\newcommand{\Pp}{\mathcal P}
\newcommand{\op}[1]{\mathsf{#1}}
\newcommand{\Dd}{\mathrm{D}}
\newcommand{\duality}[2]{\langle #1,#2\rangle}
\newcommand{\ad}{\operatorname{ad}}
\newcommand{\ord}{\operatorname{ord}}
\newcommand{\spanC}{\operatorname{span}_{\mathbb C}}

\setlist{itemsep=7pt,topsep=7pt,parsep=3pt}
\allowdisplaybreaks[1]
\title[Finite recurrences]{Polynomial eigenfunctions and recurrences of fixed length}
\author{K. Castillo}
\address{CMUC, Department of Mathematics, University of Coimbra,
3000-143 Coimbra, Portugal}
\email{math@keniercastillo.com}
\date{}
\hypersetup{pdftitle={Polynomial eigenfunctions and recurrences of fixed length},
pdfauthor={K. Castillo}}
\begin{document}
\begin{abstract}
This paper solves the fixed-length recurrence problem posed by Horozov,
Shapiro and Tater. For the rational Darboux transformations specified here,
it also resolves the order and classification questions raised by their
two conjectures. It classifies monic polynomial sequences whose members
are eigenpolynomials of one fixed nonscalar differential operator of
arbitrary finite order with polynomial coefficients and for which
multiplication by the variable involves a fixed number of consecutive
polynomials. Orthogonality is not assumed, and only sufficiently large
degrees need be prescribed. If the lowest recurrence index is $n-m+1$,
with $m\ge3$, the differential order equals $m$ times the degree of the
eigenvalue as a polynomial in $n$. After translating the variable,
every such sequence is a finite rational Darboux transform of an explicit
family with differential order $m$ and eigenvalue $n$.
Necessary and sufficient conditions on the factors determine all polynomial
differential eigenoperators of the resulting family. For families that cannot be obtained by these
transformations from families of smaller minimum differential order,
the minimum differential and recurrence orders coincide. All such families
are classified in minimum differential orders at least three, including
additional families in composite orders.
\end{abstract}
\maketitle
\begingroup
\setlength{\parskip}{0pt}
\tableofcontents
\endgroup

\section{The fixed-length recurrence problem}
\label{sec:introduction}

The problem is to classify polynomial sequences with two properties:
multiplication by the variable involves a fixed number of consecutive
polynomials, and one fixed differential operator with polynomial
coefficients acts diagonally on the sequence. The differential order
is arbitrary. The sequence may be prescribed only in sufficiently
large degrees, and orthogonality is not assumed. This is Problem~2
of Horozov, Shapiro and Tater \cite{hst2025}.
The central ideas are the same as in the work of Castillo and
Shapiro \cite[Sections~2.1--2.2, 3.3 and~4]{castilloshapiro2026}:
duality on polynomial space, evaluation and derivative functionals
to control poles, and bispectral Darboux factorizations.
Here they are applied to recurrences of arbitrary fixed length.

\subsection{The problem and the conjectures}

Write $\Pp=\mathbb C[x]$, $\partial=d/dx$, $\op E=x\partial$ and
$\op I$ for the identity. For
$\op A=\sum_{j=0}^r a_j(x)\partial^j$, $a_j\in\mathbb C(x)$, $a_r\ne0$,
its order and principal coefficient are $\ord\op A=r$ and
$\operatorname{lc}(\op A)=a_r(x)$. Scalars are $c\op I$, $c\in\mathbb C$.
Set $\Pp[\partial]=\mathbb C[x][\partial]$, $\mathcal D=\mathbb C(x)[\partial]$,
and let $\mathcal Q$ be its division ring. Order and principal coefficient extend by
$$
 \ord(\op A\op B^{-1})=\ord\op A-\ord\op B,\quad
 \operatorname{lc}(\op A\op B^{-1})=
 \frac{\operatorname{lc}(\op A)}{\operatorname{lc}(\op B)}.
$$

Let $\op L$ be a fixed nonscalar polynomial differential operator, and let
$P(x)=(P_n(x))_{n\ge n_0}$, $n_0\ge0$, be monic, with $\deg P_n(x)=n$ and
$\op L P_n(x)=\lambda_nP_n(x)$ for $n\ge n_0$, where $\lambda_n\in\mathbb C$.
Lemma~\ref{prop:unrestricted-eigenpolynomials} gives a nonconstant
polynomial $\lambda$ with $\lambda_n=\lambda(n)$.
The required recurrence is
\refstepcounter{equation}\label{eq:finite-recurrence}
$$
 xP_n(x)=P_{n+1}(x)+\sum_{s=0}^{w}b_s(n)P_{n-s}(x),\quad n\ge N,
 \leqno\textup{(\theequation)}
$$
with integers $w\ge0$, $N\ge n_0+w$ and coefficients independent of $x$.
The recurrence order is the least $w+1$ for which
\eqref{eq:finite-recurrence} holds for all sufficiently large $n$.
The extreme shifts are $1,-w$ if $w\ge1$, and the order is one if
$w=0$. Initial coefficients may vanish.

The minimum differential order is the least order of a nonscalar
polynomial differential operator having all prescribed polynomials as
eigenfunctions \cite[Remark~1.5]{hst2025}. Unless the full sequence is specified, only sufficiently
large degrees are prescribed \cite[p.~2, after~(1.2)]{hst2025}.
The differential eigenoperator algebra consists of all polynomial
differential operators having every prescribed $P_n(x)$ as an eigenfunction.
Products act from right to left.

For a sequence $v=(v_n)_{n\ge n_0}$, put
$(\op T v)_n=v_{n+1}$, $(\op T^{-1}v)_n=v_{n-1}$,
and $(rv)_n=r(n)v_n$ for $r\in\mathbb C(n)$; thus
$\op T r(n)=r(n+1)\op T$.
Finite rational difference operators are sums
$\sum_{i=p}^q r_i(n)\op T^i$, $p,q\in\mathbb Z$, $r_i\in\mathbb C(n)$.
For a nonzero operator, choose $r_p,r_q\ne0$; its width is $q-p$.
Identities are considered only at indices where all terms are defined. The recurrence operator is
$$
 \op J=\op T+\sum_{s=0}^{w}b_s(n)\op T^{-s}.
$$
The recurrence is the difference eigenvalue equation $\op J P(x)=xP(x)$
for the prescribed column $P(x)$. Together with $\op L P_n(x)=\lambda_nP_n(x)$,
this is the bispectral property used here: the operators act in
$x$ and $n$, and their eigenvalues depend on the other variable.

For a linearly independent rational family
$\psi(x)=(\psi_n(x))_{n\ge n_0}$, let $(\mathbf a_j)_{j\ge n_0}$ be its coordinate
functionals on its span: $\duality{\mathbf a_j}{\psi_n(x)}=\delta_{jn}$,
where $\delta_{jn}$ is the Kronecker delta and the pairing denotes evaluation.
Its Fourier algebra consists of $\op A\in\mathcal D$ such that
$$
 \op A\psi_n(x)=\sum_{i=p}^q r_i(n)\psi_{n+i}(x),\quad
 r_i(n)\in\mathbb C(n),
$$
for all sufficiently large $n$, with fixed $p,q$.
The Fourier map $\op b$ sends each finite difference operator admitting a
differential partner to that partner. Its inverse is
$$
 r_i(n)=\duality{\mathbf a_{n+i}}{\op A\psi_n(x)},\quad
 \op b^{-1}(\op A)=\sum_{i=p}^q r_i(n)\op T^i.
$$
Thus $\op A\psi(x)=\op b^{-1}(\op A)\psi(x)$, with coordinates in the current
family's span, including rational intermediate families.
The coordinate functionals determine $\op b^{-1}(\op A)$ uniquely. A finite rational band with
finite-dimensional image on a tail is zero: all high coordinate
functionals annihilate that image, so its rational coefficients vanish
on all large integers. The differential action is faithful because a
nonzero differential operator has finite-dimensional kernel.
The inverse Fourier map reverses order:
$$
 \op b^{-1}(\op A\op B)=\op b^{-1}(\op B)\op b^{-1}(\op A).
$$
The following shifted finite-band form of the triangular Darboux
reversal in \cite[Remark~1.6]{hst2025} is used:
$$
 \op J-a\op I=\op U\op F,\quad
 \widehat{\op J}-a\op I=\op F\op U,\quad
 \widehat\psi(x)=\op F\psi(x),
$$
where $a\in\mathbb C$; $\op U,\op F$ are finite rational difference
operators with, respectively, nonnegative and nonpositive shifts, and
$\op F$ has nonzero diagonal. Monic normalization preserves degree labels.
Every polynomial family in an allowed chain must have a fixed
nonscalar polynomial differential eigenoperator. Bispectrality follows
from the factorizations below; for recurrence order two, see also
\cite[Theorem~4.1(i) and its proof, pp.~29, 33--35]{castilloshapiro2026}.

On the differential side \cite{bhy1996},
a spectral factorization $\op C=\op Q\theta(x)^{-1}\op P$, where
$\theta(x)$ is a nonzero rational function,
$\op b^{-1}(\op C)=c(n)$ is a nonconstant rational function, and
$\op P,\op Q$ and the multiplier $\theta(x)$ belong to the Fourier algebra,
gives the family $\op P\psi(x)$ and operator $\op P\op Q\theta(x)^{-1}$.
These conditions are verified in
Proposition~\ref{prop:general-bispectral-reduction} and
Lemma~\ref{lem:listed-source-exact-bhy}.

The allowed recurrence-side bispectral Darboux steps are elementary
reversals and spectral reversals
$$
 \begin{gathered}
 \theta(\op J)=\op B\op A,\quad
 \theta(\widehat{\op J})=\op A\op B,\\[7pt]
 \widehat{\op J}\op A=\op A\op J,\quad
 \widehat\psi(x)=\op A\psi(x).
 \end{gathered}
$$
Here $\theta\in\mathbb C[z]$ is nonconstant and $\op A,\op B$ are
finite rational difference operators. The forward factor is
$\op A=\alpha(n)^{-1}\op b^{-1}(\op P)$, with $0\ne\alpha\in\mathbb C(n)$
and $\op P$ the forward differential factor of a spectral factorization
satisfying the conditions above. A spectral reversal counts as one step.
The initial family in a chain is called its source.
A solution is Darboux irreducible if no finite composition of these
steps produces it from a polynomial eigenfamily with a finite
recurrence for $x$ and smaller minimum differential order. Both minima
refer to the prescribed degrees.

Conjecture~1.9 of \cite{hst2025} asserts equality of the minimum
differential and recurrence orders for irreducible solutions.
The corresponding assertion is proved here for the steps defined
above, including spectral reversals with nontriangular factors. The assertion for chains
restricted to the triangular reversals of \cite[Remark~1.6]{hst2025}
requires a separate argument.

The families proposed for irreducible solutions of minimum order
$r>0$ in \cite[Conjecture~1.10]{hst2025} have, up to invertible affine
changes of $x$ and the eigenvalue, the forms
$$
 \op E+\sum_{j=1}^r a_jx^{j-1}\partial^j,\quad a_r\ne0,
$$
and
$$
 \op E+Q'(\op G)\op G,\quad
 \op G=R(\op E)\partial,\quad
 \deg R=\ell-1,\quad \deg Q=r/\ell,\quad Q(0)=0.
$$
Here $a_j\in\mathbb C$, $Q,R\in\mathbb C[z]$, and $\ell$ is a
positive divisor of $r$. The first type is contained in the second.
The operator $\op G$ is the generator; it is called polynomial when
its coefficients are polynomial. In the classification below, $\op G$
may have rational coefficients, provided that each power occurring in
$\op L$ has polynomial coefficients.

\subsection{Main results}

Recurrence orders one and two follow from
\cite[Theorem~4.1, p.~29]{castilloshapiro2026}.

\begin{theorem}
\label{thm:all-polynomial-sources}
Let $\op L$ be a nonscalar polynomial differential operator, and let
$P(x)=(P_n(x))_{n\ge n_0}$ be monic, with $\deg P_n(x)=n$.
Suppose that its members are eigenpolynomials of $\op L$ in all
sufficiently large degrees and that its recurrence order is $m\ge3$.
In these degrees, $\op L P_n(x)=\lambda(n)P_n(x)$ for a nonconstant polynomial $\lambda$, and
$\ord\op L=m\deg\lambda$.
After a translation of $x$, the prescribed family is a finite rational
bispectral Darboux transform of the full monic eigenfamily
$\psi(x)=(\psi_n(x))_{n\ge0}$ of
$$
 \op K=\op E+q(\op G),\quad
 \op G=x^{-1}g(\op E),\quad
 g(0)=q(0)=0,\quad
 (\deg q)(\deg g)=m,
$$
where $q,g$ are nonconstant polynomials and $g$ has no positive integer root.
The operator $\op K$ has polynomial coefficients and order $m$, with
$\op K\psi_n(x)=n\psi_n(x)$.
The family $\psi(x)$ has recurrence order $m$, and its last
recurrence coefficient is nonzero for every $n\ge m-1$.
The transformation sends $\psi_n(x)$ to a nonzero multiple of $P_n(x)$
at the same index $n$, for all sufficiently large $n$.
Its forward differential operator intertwines $\lambda(\op K)$
with $\op L$. The source and Darboux operators are fixed.
\end{theorem}

\begin{theorem}
\label{thm:general-irreducible-order-equality}
Let $(P_n(x))_{n\ge n_0}$ be monic, with $\deg P_n(x)=n$.
Suppose that in all sufficiently large degrees it has recurrence order
$m$ and a nonscalar polynomial differential eigenoperator. Let $r$ be
the minimum differential order on these prescribed degrees.
If this prescribed tail is Darboux irreducible, then $r=m$.
\end{theorem}

A rational generator gives the polynomial operator
\refstepcounter{equation}\label{eq:singular-fourth-order-reformulation}
$$
 \op L=\op E+\bigl(x^{-1}(\op E^2-1)\bigr)^2
       =\op E+x^2\partial^4+4x\partial^3,
 \leqno\textup{(\theequation)}
$$
whose monic eigenfamily is Darboux irreducible and has minimum
differential and recurrence orders four. Section~\ref{sec:principal-low-order} gives
families of this form in every composite order.

\begin{theorem}
\label{thm:general-irreducible-classification}
Let $(P_n(x))_{n\ge n_0}$ be monic, with $\deg P_n(x)=n$,
and suppose that it has a recurrence of fixed length in sufficiently
large degrees. Let $\op L$ be a polynomial differential eigenoperator
of minimum order $r\ge3$ on these prescribed degrees.
The prescribed tail is Darboux irreducible exactly when, after invertible affine
changes of $x$ and of the eigenvalue, $\op L$ has the form
$$
 \op L=\op E+\sum_{i\in S}c_i\op G^i,\quad
 \op G=x^{-1}g(\op E),\quad
 K\deg g=r,\quad K=\max S,
$$
where $S$ is a finite nonempty set of positive integers,
$c_i\ne0$, $g(z)\in\mathbb C[z]$ is nonconstant, and each
occurring power $\op G^i$ has polynomial coefficients.

After these affine changes and monic normalization, the prescribed
family is a tail of the full monic eigenfamily
$$
 P_n(x)=\exp\left(\sum_{i\in S}\frac{c_i}{i}\op G^i\right)x^n,
 \quad n\ge0.
$$
The exponential is finite on every polynomial, and the recurrence
order is $r$.
\end{theorem}

Removing the common monomial factor and shifting degree gives a
polynomial generator of the same minimum order
(Proposition~\ref{prop:regular-core-normal-form}). A polynomial generator
already suffices if the full normalized family's last recurrence
coefficient is nonzero for every $n\ge r-1$.

For polynomial differential operators $\op K$ and $\op F\ne0$, set
$$
 \mathcal A_F=\{\lambda\in\mathbb C[z]:
 \op F\lambda(\op K)\op F^{-1}\in\Pp[\partial]\}.
$$

\begin{theorem}
\label{thm:complete-admissible-data}
Fix $m\ge3$. Choose nonconstant polynomials $g(z)$ and
$q(z)=\sum_{j=1}^k q_jz^j$, where $k=\deg q$, with
$$
 q(0)=g(0)=0,\quad (\deg q)(\deg g)=m.
$$
Positive integer roots of $g$ are allowed. Put
$$
 \begin{gathered}
 \op G=x^{-1}g(\op E),\quad \op K=\op E+q(\op G),\quad
 \beta=\max\{j\in\mathbb Z_{\ge0}:g(-j)=0\},\\[7pt]
 \psi_n(x)=\exp\left(\sum_{j=1}^k\frac{q_j}{j}\op G^j\right)x^n,
 \quad n\ge0.
 \end{gathered}
$$
The exponential is finite on every monomial, and
$\op K\psi_n(x)=n\psi_n(x)$ for $n\ge0$.
Choose $s\in\mathbb Z_{\ge0}$ and polynomials $t_j(z)$,
$0\le j\le s$, with $t_0\ne0$, and put
$$
 \op F=\sum_{j=0}^s\op G^j t_j(\op E),\quad \alpha(n)=t_0(n).
$$
For sufficiently large $N$, put
$$
 P(x)=(P_n(x))_{n\ge N},\quad
 P_n(x)=\alpha(n)^{-1}\op F\psi_n(x).
$$
This sequence is monic, with $\deg P_n(x)=n$.
Let $\op b$ be the Fourier map of $\psi(x)=(\psi_n(x))_{n\ge0}$ and put
$$
 \begin{aligned}
 \op B&=g(n)\op T^{-1},\quad \op b^{-1}(\op G)=\op B,\quad
 \op b^{-1}(\op E)=n-q(\op B),\\[7pt]
 \op J&=\op b^{-1}(x)=g(n+1)^{-1}\op T g\bigl(n-q(\op B)\bigr),\\[7pt]
 \op A&=\alpha(n)^{-1}\op b^{-1}(\op F).
 \end{aligned}
$$
The family $P(x)$ has a nonscalar polynomial differential
eigenoperator and a recurrence of fixed length for multiplication by $x$
exactly when both conditions hold:
\begin{enumerate}[label=(\roman*)]
\item There are a nonconstant polynomial $c(z)$ and an operator
$\op S\in\mathbb C[x][\partial]$ such that
\refstepcounter{equation}\label{eq:admissible-data-factorization}
$$
 x^\beta c(\op K)=\op S\op F.
 \leqno\textup{(\theequation)}
$$
\item There is a finite rational difference operator $\widehat{\op J}$ such that
$$
 \widehat{\op J}\op A=\op A\op J.
$$
\end{enumerate}
Under these conditions, $xP(x)=\widehat{\op J}P(x)$ has recurrence order $m$,
where $P(x)$ is the prescribed column.
The family is a finite rational bispectral Darboux transform of $\psi(x)$.
Up to translation of $x$, every solution of recurrence order $m$ in
sufficiently large degrees is obtained in this way.

For each choice satisfying \textup{(i)} and \textup{(ii)}, all polynomial
differential eigenoperators for the prescribed degrees are
$$
 \op L_\lambda=\op F\lambda(\op K)\op F^{-1},\quad
 \lambda\in\mathcal A_F.
$$
If $\lambda$ is nonconstant, then $\ord\op L_\lambda=m\deg\lambda$.
The finite system in Corollary~\ref{cor:finite-spectral-jet-algebra}
determines $\mathcal A_F$. Lemma~\ref{lem:admissible-recurrence-test}
gives the coefficient test for \textup{(ii)}. The minimum differential order is $m$ times
the least degree of a nonconstant polynomial in this algebra.
The greatest common divisor of the positive differential orders is
$m$, and every sufficiently large multiple of $m$ occurs.
\end{theorem}

\subsection{Earlier results}

The Bochner--Krall problem, Problem~1 of Horozov, Shapiro and Tater, is treated in
\cite[Theorem~4.1, p.~29]{castilloshapiro2026}, which covers
recurrence orders one and two. Krall's problem \cite[pp.~705--706]{krall1938} requires
orthogonal polynomials with nonzero squared norms,
hence a three-term recurrence. Its classical sources are Hermite,
Laguerre, Jacobi and Bessel; the nonclassical families are prescribed
Darboux transforms of Laguerre or Jacobi
\cite[Theorem~2.1(ii), p.~7]{castilloshapiro2026}.

The proofs use Maroni's functional approach
\cite{maroni1985,maroni1988,maroni1991a}, as developed for the
three-term problem in \cite[Sections~2.1--2.2 and~3.3]{castilloshapiro2026}.
Transposition on the coordinate functionals proves the order bound
and order equality. Evaluation and derivative functionals give the
conditions on the Darboux factors. The source classification uses
operators commuting with homogeneous operators \cite{bavula2005},
Horozov's constructions \cite{horozov2016,horozov2018vector},
and local expansions. Shapiro's book develops the basic theory of
exactly solvable operators and presents the generalized Bochner--Krall
problem and its conjectures \cite[Sections~1.2--1.3 and~5.3]{shapiro2026book}.

Write $\mathbb C[[z]]$ for the ring of formal power series and
$\mathbb C((z))$ for its field of formal Laurent series.
An Appell family satisfies $\partial P_n(x)=nP_{n-1}(x)$ for $n\ge1$.
A monic Sheffer family has exponential generating function
$A(z)e^{x\psi(z)}$, where $A,\psi\in\mathbb C[[z]]$,
$A(0)=\psi'(0)=1$ and $\psi(0)=0$.

Sheffer proved the Appell existence and minimum-order criteria
\cite[Theorems~4--5]{sheffer1935appell}. He also gave a canonical form
and a finite-order criterion for Sheffer eigenoperators
\cite[Theorems~3.4, 3.6 and Corollary~3.3]{sheffer1939}.
Krall and Sheffer treated orthogonal polynomial eigenfamilies for possibly
infinite-order equations \cite{krallsheffer1966}.
Horozov's constructions are fundamental here: they provide the polynomial
sources, their differential and recurrence orders, and their raising and
lowering operators \cite[Theorems~4.1--4.2 and Remark~4.3]{horozov2018vector}.
The converse uses Bavula's results on operators commuting with a fixed
homogeneous operator \cite{bavula2005}.

Horozov, Shapiro and Tater classified irreducible operators of
differential and recurrence orders three
\cite[Theorem~1.14]{hst2025}. For differential order three, the
extension to arbitrary recurrence order was posed in
\cite[Final remarks, II]{hst2025}. For the Darboux
steps used here, Theorems~\ref{thm:general-irreducible-order-equality}
and~\ref{thm:general-irreducible-classification} give the classification
in every minimum differential order at least three. Related classifications of third-order eigenfamilies appear in
\cite[Propositions~2 and~4]{mesquita2022} and
\cite[Theorem~3]{mesquita2024}. Theorem~\ref{prop:homogeneous-recurrence} extends the
monomial classification of \cite{abst2026} to arbitrary polynomial
eigenvalues and homogeneous lowering terms. For zero asymptotics and algebraic
curves of that monomial family, see \cite{shapirotater}.

Casper and Yakimov prove a matrix Bochner theorem
\cite[Theorem~1.3]{casperyakimov2022} for a class of matrix-valued
orthogonal polynomials with a second-order eigenoperator symmetric
for the matrix weight. Their hypotheses include positive definiteness
of its principal coefficient times the weight. Here the families are
scalar and the differential order is arbitrary.

Bakalov, Horozov and Yakimov established a systematic Darboux framework
for constructing bispectral operators \cite{bhy1996,bhy1996grassmann}
and bispectral algebras of arbitrary rank \cite{bhy1997}.
Their algorithmic construction unifies and extends the earlier approaches
of Duistermaat and Gr\"unbaum \cite{duistermaatgrunbaum1986}
and Wilson \cite{wilson1993}. Horozov subsequently obtained the complete
classification of normalized bispectral differential operators of prime
order \cite{horozov2002prime}, recovering the second-order classification
as a special case. Rank is the greatest common divisor of differential
orders. Here the spectral variable is the integer degree and the
recurrence represents multiplication by $x$.

\section{Polynomial space and duality}
\label{sec:preliminaries}

The functional definitions and weak dual expansions are those of
\cite[Sections~2.1--2.2, pp.~6--7]{castilloshapiro2026}, with coordinate
functionals for an arbitrary monic basis. Put
$\Pp_k=\{p\in\Pp:\deg p\le k\}$ for $k\ge0$ and $\deg0=-\infty$.
Write
$$
 [\op A,\op B]=\op A\op B-\op B\op A,\quad
 \ad_{\op A}(\op B)=[\op A,\op B],\quad
 z^{\underline j}=\prod_{i=0}^{j-1}(z-i),\quad z^{\underline0}=1.
$$
Write $[z^j]f(z)$ for the coefficient of $z^j$.

Give $\Pp=\bigcup_k\Pp_k$ its locally convex inductive-limit topology.
Every linear functional and endomorphism is continuous. Give its dual
$\Pp'$ the weak topology $\sigma(\Pp',\Pp)$, with pairing
$\duality{\mathbf u}{p(x)}=\mathbf u(p(x))$. A monic basis $(P_n(x))_{n\ge0}$ has a dual sequence
$(\mathbf a_n)_{n\ge0}$ such that, for every $\boldsymbol\eta\in\Pp'$,
$$
 \duality{\mathbf a_n}{P_j(x)}=\delta_{nj},\quad
 \boldsymbol\eta=\sum_{n\ge0}
 \duality{\boldsymbol\eta}{P_n(x)}\mathbf a_n.
$$
The expansion converges weakly, since $\duality{\mathbf a_n}{p(x)}=0$
for $n>\deg p$. Thus $\Pp'$ is the product of its scalar coordinates
with their usual complex topologies. Complete omitted degrees by arbitrary
monic polynomials before taking the dual. Coefficient localizations below
are algebraic; no infinite dual expansion passes through a localization.
For $f,p\in\Pp$, $\mathbf u\in\Pp'$ and linear maps
$\op A,\op B:\Pp\to\Pp$, define
$$
 \begin{aligned}
 \duality{f(x)\mathbf u}{p(x)}&=\duality{\mathbf u}{f(x)p(x)},\\[7pt]
 \duality{\Dd\mathbf u}{p(x)}&=-\duality{\mathbf u}{\partial p(x)},\\[7pt]
 \duality{\op A'\mathbf u}{p(x)}&=\duality{\mathbf u}{\op A p(x)}.
 \end{aligned}
$$
Then $\Dd=-\partial'$, $(\op A\op B)'=\op B'\op A'$, and
$$
 \Dd(f(x)\mathbf u)=f(x)\Dd\mathbf u+(\partial f(x))\mathbf u.
$$
Every transpose is weakly continuous. Pairings are bilinear, without
conjugation. For finite-dimensional complex spaces, write
$V^\vee=\operatorname{Hom}_{\mathbb C}(V,\mathbb C)$ and give it
its usual topology. The transpose of $A:V\to W$ is
$A^\vee:W^\vee\to V^\vee$, defined by
$\langle A^\vee\ell,v\rangle=\langle\ell,Av\rangle$.
For $\boldsymbol\delta_a(p(x))=p(a)$,
$\duality{\Dd^j\boldsymbol\delta_a}{p(x)}=(-1)^jp^{(j)}(a)$.
For these topologies and dual expansions, see
\cite[Sections~1--2]{maroni1985}, \cite[Section~1]{maroni1988} and
\cite[Theorem~1.3 and Section~A.6]{castillopetronilho2ed}.
For recurrences holding in every degree with $w\ge1$ and
$b_w(n)\ne0$ for every $n\ge w$, see
\cite[Theorem~2.1]{maroni1989}. For the vector Pad\'e formulation
under additional initial nonvanishing conditions, see
\cite[Chapter~II, Theorems~III.2--III.3]{vaniseghem1987}.
No orthogonality functional is needed.
Power series are compared coefficient by coefficient.
A rational differential operator $\op A=\sum_{j=0}^s a_j(x)\partial^j$
mapping $\Pp$ into $\Pp$ has polynomial coefficients: evaluation on
$1,x,\ldots,x^s$ proves this successively for $a_0,\ldots,a_s$
\cite[proof of Lemma~3.6, p.~21]{castilloshapiro2026}.
A finite point is ordinary when all coefficients are holomorphic and the
principal coefficient is nonzero. An order-$r$ homogeneous equation,
$r\ge1$, has an $r$-dimensional local solution space there, parametrized
by its first $r$ initial values.

Let $\op A,\op T\in\mathcal D$, $\op T\ne0$, and suppose that
$V=\ker\op T$ at a common ordinary point is $\op A$-invariant.
For $f\in\mathbb C[z]$, the identity
$f((\op A|_V)^\vee)=0$ implies
$$
 f(\op A)\in\mathcal D\op T.
$$
Indeed, dual separation gives $f(\op A)V=0$, and a nonzero
right-division remainder has order below $\dim V$ and cannot annihilate $V$.
Cayley--Hamilton supplies such a nonzero $f$.

\begin{lemma}
\label{prop:unrestricted-eigenpolynomials}
Let $\op L=\sum_{j=0}^r a_j(x)\partial^j$ be a nonscalar polynomial
differential operator of positive order. It has a monic eigenpolynomial
of every sufficiently large degree exactly when $\deg a_j\le j$ for
all $j$ and
$$
 \lambda(z)=\sum_{j=0}^r[x^j]a_j(x)z^{\underline j}
$$
is nonconstant. These monic eigenpolynomials are unique in every
sufficiently large degree. For each fixed $k\ge0$, the coefficient
of $x^{n-k}$ is a rational function of $n$.
\end{lemma}

\begin{proof}
Put $t=\max_j(\deg a_j-j)$. If $t>0$, then
$$
 [x^{n+t}]\op L P_n(x)
 =\sum_{\deg a_j-j=t}[x^{j+t}]a_j(x) n^{\underline j}.
$$
This nonzero polynomial cannot vanish for all large $n$. Hence
$\deg a_j\le j$, and comparison at degree $n$ gives $\lambda(n)$.
If $\lambda$ were constant, $\op L-\lambda\op I$ would be nonzero
and annihilate infinitely many independent polynomials, a contradiction.

Conversely, $\op L$ is triangular on monomials, with diagonal
$\lambda(n)$. The real polynomial $|\lambda(n)|^2$ has positive leading
coefficient and is eventually strictly increasing. Thus
$\lambda(n)\ne\lambda(j)$ for all large $n$ and $0\le j<n$.
Descending coefficient comparison gives a unique monic eigenpolynomial.
For fixed $k$, its coefficient of $x^{n-k}$ uses only polynomial
expressions in $n$ and division by $\lambda(n)-\lambda(n-i)$,
$1\le i\le k$, so it is rational. For the triangular comparison, see
\cite[proof of Lemma~3.1, p.~9]{castilloshapiro2026} and
\cite[Section~1.2, pp.~3--4]{shapiro2026book}. Eventual existence
and uniqueness are also stated in \cite[p.~2, following~(1.2)]{hst2025}.
\end{proof}

Triangular comparison in \eqref{eq:finite-recurrence} shows that each
$b_s(n)$ agrees with a rational function for all sufficiently large $n$.
Nonzero rational functions have finitely many zeros and poles, so the
extreme shifts are fixed beyond a finite initial block.

\begin{theorem}
\label{thm:general-order-bound}
Let $\op{L}$ be a nonscalar polynomial differential operator of order
$r$, with monic eigenpolynomials in every sufficiently large degree.
If these polynomials satisfy \eqref{eq:finite-recurrence}, then the
recurrence order is at most $r$.
\end{theorem}

\begin{proof}
The case $w=0$ is immediate, so suppose $w\ge1$ and $b_w\ne0$.
Complete the omitted degrees and take the dual basis $(\mathbf a_n)_{n\ge0}$.
For $\op L=\sum_{j=0}^r a_j(x)\partial^j$,
$\op L'=\sum_{j=0}^r(-\Dd)^j a_j(x)$, and in sufficiently large indices,
$$
 x\mathbf a_n=\mathbf a_{n-1}
 +\sum_{s=0}^w b_s(n+s)\mathbf a_{n+s},\quad
 \op L'\mathbf a_n=\lambda(n)\mathbf a_n.
$$
Every $x\mathbf a_n$ is a finite dual sum. Choose $H$ beyond all
exceptional indices and supports. For $0\le j<w$ and $q\ge0$,
$x^q\mathbf a_{H+j}$ has nonzero coefficient at its largest index
$H+j+qw$. These indices are distinct, so
$\mathbf a_H,\ldots,\mathbf a_{H+w-1}$ are independent over
$\mathbb C[x]$. Localize algebraically:
$$
 \mathcal M=\mathbb C(x)\otimes_{\mathbb C[x]}\Pp',
$$
where $\mathbf u$ denotes its image $1\otimes\mathbf u$.
They remain independent over $\mathbb C(x)$: clearing denominators
and multiplying by a nonzero polynomial would give a polynomial relation.
This holds for every large $H$. For $r(x)\in\mathbb C(x)$, extend $\Dd$ by
$$
 \Dd(r(x)\otimes\mathbf u)
 =(\partial r(x))\otimes\mathbf u+r(x)\otimes\Dd\mathbf u.
$$

For a nonzero polynomial $f$, put
$\op B_f=f\op L'-(\op Lf)\op I$. For every $\mathbf u\in\Pp'$,
$$
 \duality{\op B_f\mathbf u}{1}
 =\duality{\mathbf u}{\op Lf-\op Lf}=0.
$$
Division with $\Dd$ on the left has a multiplication remainder;
the same pairing forces it to vanish. Thus
$$
 \op B_f=\Dd\op C_f,\quad \ord\op C_f=r-1,
 \quad\operatorname{lc}(\op C_f)=(-1)^rf a_r.
$$ For $m=H+w-1$ and $H$
sufficiently large, $\Dd\op C_{P_m}\mathbf a_m=0$ in $\Pp'$.
Surjectivity of $\partial:\Pp\to\Pp$ makes $\Dd$ injective there, so
$\op C_{P_m}\mathbf a_m=0$. Hence
$$
 \mathcal W=\mathbb C(x)[\Dd]\mathbf a_m\subseteq\mathcal M
$$
has dimension at most $r-1$ over $\mathbb C(x)$ and is stable under
$\Dd$, rational multiplication, $\op L'$ and $x$. By the preceding lemma,
the eigenvalues on each sufficiently high finite range of indices are distinct.
Interpolation in $\op L'$ on $x\mathbf a_n$ therefore extracts
$\mathbf a_{n-1}$ when $\mathbf a_n\in\mathcal W$ and $n$ is large.
Descending from $m$ gives $\mathbf a_H,\ldots,\mathbf a_{H+w-1}\in\mathcal W$,
so $w\le\dim_{\mathbb C(x)}\mathcal W\le r-1$.
\end{proof}

\begin{corollary}
\label{cor:width-irreducibility}
Let a polynomial eigenfamily have differential and recurrence orders
equal to $N\ge2$. It cannot arise from a polynomial differential
source of order below $N$ by a finite rational bispectral Darboux
transformation between solutions of the fixed-length problem.
\end{corollary}

\begin{proof}
The forward Fourier factors and monic normalizations give a nonzero
finite rational difference operator $\op A$ with $\widehat\psi(x)=\op A\psi(x)$.
Multiplication by $x$ and uniqueness of finite difference partners give
$\op J\op A=\op A\op J_0$, also for factorizations of $\theta(\op J_0)$.
Extreme shifts add in products, each from a unique nonzero pair.
Hence $\op J$ and $\op J_0$ have width $N$, and
Theorem~\ref{thm:general-order-bound} implies that the source has order at least $N$.
\end{proof}

\begin{lemma}
\label{lem:order-one-canonical}
If a monic polynomial eigenfamily of a nonscalar polynomial
differential operator has recurrence order one, then in all
sufficiently large degrees
$P_n(x)=(x-\beta)^n$ for some $\beta\in\mathbb C$. The minimum
nonscalar differential order for these sufficiently large degrees is one.
\end{lemma}

\begin{proof}
Apply \cite[Theorem~4.1(ii), p.~29]{castilloshapiro2026} to
$xP_n(x)=P_{n+1}(x)+b_0(n)P_n(x)$.
\end{proof}

\begin{theorem}
\label{thm:spectral-degree-divisibility}
Let a polynomial differential operator of order $r$ have monic
eigenpolynomials in every sufficiently large degree, with eigenvalue
polynomial $\lambda$ of degree $d\ge1$. If their recurrence has
fixed order different from two, then $d$ divides $r$. Moreover, its
principal coefficient is $\alpha a(x)^d$, where
$\alpha=[z^d]\lambda(z)$ and $a(x)$ is polynomial.
The same conclusions hold for a genuine three-term recurrence if
no $\zeta\in\mathbb C$ makes $\lambda(\zeta+z)$ even.
\end{theorem}

\begin{proof}
For recurrence order one, Lemma~\ref{lem:order-one-canonical} and
\cite[proof of Theorem~4.1(ii), pp.~31--32]{castilloshapiro2026} give
$\op L=\lambda((x-\beta)\partial)$ for some $\beta\in\mathbb C$.
Thus $r=d$ and the principal coefficient is $\alpha(x-\beta)^d$.
Assume henceforth that the recurrence order is at least two.

For $h\in\mathbb Z\setminus\{0\}$, let $C_h\subset\mathbb C^2$
be the algebraic curve parametrized by $(\lambda(t),\lambda(t+h))$.
This parametrization has a rational inverse: for $d>1$, the degree
$[\mathbb C(t):\mathbb C(\lambda(t),\lambda(t+h))]$ divides both $d$
and, using the difference of the coordinates, $d-1$;
the case $d=1$ is immediate.
If $C_h=C_k$, the parameter change preserves $\lambda$ and fixes its
unique pole, infinity, so it is affine. After translating the parameter to remove the term of degree $d-1$
in $\lambda$, it is
$t\mapsto\omega t$, with $\omega^d=1$. Comparison at degree $d-1$ gives
$k=\omega h$. Integer shifts therefore satisfy $h=k$ or $h=-k$,
the latter requiring reflection symmetry. The diagonal curve differs
from every $C_h$.

If the lowest recurrence shift is $-w$ with $w\ge2$, take $h=-w$;
its partner $+w$ is absent because the largest shift is one.
For $w=1$ without reflection symmetry, take $h=1$.
Let $F(X,Y)$ be the product of defining polynomials of the
curves of the other terms present in the recurrence. On polynomial
differential operators set
$\mathcal R(\op A)=\op A\op L$ and
$\mathcal S(\op A)=\op L\op A$. These actions commute, and
$\op C=F(\mathcal R,\mathcal S)(x)$ is nonzero, with
$$
 \op C P_n(x)=\rho(n)P_{n+h}(x),\quad \rho(n)\ne0
 \quad\text{for all sufficiently large }n.
$$

Let $\mathcal F$ be the ideal of finite-rank endomorphisms of $\Pp$.
In $\operatorname{End}_{\mathbb C}(\Pp)/\mathcal F$, a diagonal class
with entries $\nu(n)$ has a representative sending
$P_n(x)$ to $\nu(n)P_n(x)$ for all sufficiently large $n$.
The division ring $\mathcal Q$ embeds in this quotient. For a nonzero polynomial differential operator,
the leading degree term gives finite kernel and cokernel, hence an
invertible class; its infinite rank makes the representation injective.
Complete the omitted degrees. An inverse of $\op C$ modulo finite rank
sends $P_n(x)$ to $\rho(n-h)^{-1}P_{n-h}(x)$ for all large $n$.
Thus $\op L$ and $\op B=\op C^{-1}\op L\op C$ have diagonal classes
with entries $\lambda(n)$ and $\lambda(n+h)$. Their commutation and
injectivity give $[\op L,\op B]=0$ in $\mathcal Q$.

Choose a rational function $R(X,Y)$ with $R(\lambda(t),\lambda(t+h))=t$. Its denominator
is a nonzero polynomial in $t$, hence invertible modulo finite rank.
Then $\op N=R(\op L,\op B)\in\mathcal Q$ has diagonal class $n$ and
$\lambda(\op N)=\op L$. Since $\ord\op L>0$,
$k=\ord\op N>0$ and $r=dk$. If $a(x)$ is the leading coefficient
of $\op N$, then $a_r(x)=\alpha a(x)^d$. Since $a_r$ is polynomial, $a$ has no finite poles, so $a$ is polynomial.
\end{proof}

For the commutator condition for differential operators of order three, see
\cite[Lemma~2.1]{hst2025} and \cite[Section~5.4]{shapiro2026book}.

\begin{lemma}
\label{lem:dual-spectral-support}
Let $\op L$ be a polynomial differential operator with monic
eigenpolynomials $P_n(x)$ and eigenvalues $\lambda(n)$ in every
sufficiently large degree. Complete the lower degrees to a monic
basis, and let $(\mathbf a_n)_{n\ge0}$ be its dual sequence. On the weakly
continuous operators on $\Pp'$ put
$$
 \Delta(\op C)=-[\op L',\op C],\quad
 \mathcal S(\op C)=\op L'\op C.
$$
These two maps commute. For every $F\in\mathbb C[z]$ and every
sufficiently large $n$ one has the weak dual expansion
$$
 F(\Delta)(x)\mathbf a_n
 =\sum_{j\ge n-1}
 F\bigl(\lambda(n)-\lambda(j)\bigr)
 \duality{\mathbf a_n}{xP_j(x)}\mathbf a_j.
$$
If a polynomial differential operator $\op C$ satisfies
$\op C'\mathbf a_n=0$ for every sufficiently large $n$, then
$\op C=0$.

Let $J$ be the finite set of shifts with nonzero recurrence coefficients,
and suppose
$d=\deg\lambda\ge1$. Define
$$
 R_j(S,Y)=\prod_{\lambda(t)=Y}
 \bigl(S-\lambda(t+j)+\lambda(t)\bigr),\quad j\in J,
$$
with roots counted with multiplicity. Then
$$
 R_j(S,Y)=S^d+\sum_{k=1}^d c_{jk}(Y)S^{d-k},\quad
 c_{jk}\in\mathbb C[Y],\quad \deg c_{jk}\le k-1,
$$
and the following identity holds on $\Pp'$:
$$
 \prod_{j\in J}R_j(\Delta,\mathcal S)(x)=0.
$$
In particular, when $\lambda(n)=n$, for each positive integer $r$
the finite recurrence of order at most $r$ is equivalent to
$$
 F_r(\Delta)(x)=0,\quad
 F_r(z)=\prod_{h=1-r}^{1}(z-h).
$$
In this affine case the recurrence holds in every degree for the
unique full monic eigenbasis, with negative-index terms omitted.
\end{lemma}

\begin{proof}
Degree preservation and the prescribed equations give
$\op L'\mathbf a_n=\lambda(n)\mathbf a_n$ for all large $n$.
Also $\duality{\mathbf a_n}{xP_j(x)}=0$ for $j<n-1$, so all indices
in the weak expansion of $x\mathbf a_n$ are large when $n$ is large.
Weak continuity and the product rule for transposes give the formula
by induction on $\deg F$. In particular,
$$
 \bigl(F(\ad_{\op L})(x)\bigr)'=F(\Delta)(x).
$$
If $\op C'$ annihilates a dual tail, then $\op C\Pp$ is
finite-dimensional. A nonzero differential operator has finite kernel
and infinite-dimensional image, so $\op C=0$.

Write the recurrence as
$xP_n(x)=\sum_{j\in J}a_j(n)P_{n+j}(x)$. Its transpose is
$$
 x\mathbf a_n=\sum_{j\in J}a_j(n-j)\mathbf a_{n-j}
$$
for all large $n$. On the $j$th term, $\Delta$ and $\mathcal S$
give $\lambda(n)-\lambda(n-j)$ and $\lambda(n-j)$, so $R_j$ annihilates
that term. The commuting product vanishes on every large $\mathbf a_n$;
being the transpose of a polynomial differential operator, it is zero
by the preceding criterion. Symmetry in the roots of $\lambda(t)-Y$
gives polynomial coefficients in $Y$. These roots are $O(Y^{1/d})$
and each spectral difference has degree at most $d-1$ in $t$, so
$\deg c_{jk}\le\lfloor k(d-1)/d\rfloor\le k-1$.

For $\lambda(n)=n$, triangularity gives the unique full monic eigenbasis,
and the first formula holds for every $n$.
Since $F_r(n-j)=0$ exactly for $n-1\le j\le n+r-1$,
uniqueness of weak dual expansions identifies $F_r(\Delta)(x)=0$
with the full recurrence. An eventual recurrence gives the same identity
by the zero-operator criterion.
\end{proof}

The following lemma extends the argument for coefficient poles in
\cite[Lemma~3.6, p.~21]{castilloshapiro2026} to recurrences of arbitrary
fixed length holding in every degree.

\begin{lemma}
\label{lem:dual-pole-polynomiality}
Let $(\psi_n(x))_{n\ge0}$ be a monic polynomial basis satisfying
\eqref{eq:finite-recurrence} in every degree, with negative-index
terms omitted, $w\ge1$, and $b_w(n)\ne0$ for every $n\ge w$.
If $\op A\in\mathbb C(x)[\partial]$ and $\op A\psi_n(x)$ is
polynomial for every sufficiently large $n$, then $\op A$ has
polynomial coefficients.
\end{lemma}

\begin{proof}
If $\mathbf u\ne0$ is a finite dual sum with largest index $h$,
then $x^d\mathbf u$, $d\ge0$, has largest index $h+dw$, with nonzero coefficient.
Thus no nonzero polynomial annihilates $\mathbf u$.
At a pole $a$ of a coefficient, each principal-part coefficient of $\op A p$
is a finite linear combination of derivatives of $p$ at $a$.
Its functional is killed by a power of $x-a$ and vanishes on all large
$\psi_n(x)$, hence is a finite dual sum \cite[Lemma~1.1]{maroni1989} and is zero.
Therefore $\op A\Pp\subset\Pp$, and the coefficient criterion applies.
\end{proof}

\section{Order equality and bispectral reduction}
\label{sec:general-reduction}

\begin{lemma}
\label{lem:ordinary-root-normalization}
Let $\op N\in\mathcal Q$ and let $\lambda(z)\in\mathbb C[z]$ have
degree $d\ge1$. Suppose that $\lambda(\op N)=\op L\in\mathcal D$
has positive order $r$. There are nonzero $\op T,\op R\in\mathcal D$,
an operator $\op H\in\mathcal D$ of order $r/d$, and a nonconstant
polynomial $F(z)\in\mathbb C[z]$ such that
$$
 \begin{gathered}
 \op T\op N=\op H\op T,\quad
 \op T\op L=\lambda(\op H)\op T,\\[7pt]
 \op R\op T=F(\op L),\quad
 \op T\op R=F(\lambda(\op H)),\quad
 \op L\op R=\op R\lambda(\op H).
 \end{gathered}
$$
\end{lemma}

\begin{proof}
Choose $\op S\ne0$ with $\op S\op N^j\in\mathcal D$ for
$0\le j<d$. Reduction by $\lambda(\op N)=\op L$ gives this for
all $j\ge0$. The nonzero left ideal
$$
 I=\sum_{j\ge0}\mathcal D\op S\op N^j
$$
is stable under right multiplication by $\op N$. Since $\mathcal D$ is
left principal, write $I=\mathcal D\op T$. Then $\op T\op N=\op H\op T$ for some
$\op H\in\mathcal D$, of order $r/d$, and
$\op T\op L=\lambda(\op H)\op T$.

At a common ordinary point, $\op L$ preserves $V=\ker\op T$.
For a nonconstant multiple $F$ of the characteristic polynomial of
$(\op L|_V)^\vee$, Section~\ref{sec:preliminaries} gives
$F(\op L)=\op R\op T$.
Multiplication by $\op T$ and cancellation give the remaining identities.
\end{proof}

For an indeterminate $\eta$, a differential module over
$\mathbb C(\eta,x)$ is a vector space
with a $\partial$-action satisfying the product rule. It is simple if
no proper nonzero subspace is stable under $\partial$.

\begin{lemma}
\label{lem:generic-eigenvalue-simple}
Let $\op H\in\mathcal D$ have positive order $k$, let $\eta$ be an
indeterminate, and put $\mathcal D_\eta=\mathbb C(\eta,x)[\partial]$.
The left module
$$
 \mathcal D_\eta/\mathcal D_\eta(\op H-\eta)
$$
is simple and has dimension $k$ over $\mathbb C(\eta,x)$.
\end{lemma}

\begin{proof}
Write $[\op A]$ for the class of $\op A$ in this quotient.
Let $W$ be a nonzero proper submodule. Clearing a representative's
denominators in $\eta,x$ gives a nonzero class in $W$ represented by
$\sum_i\op A_i\eta^i$, with $\op A_i\in\mathcal D$. Its class equals
that of $\sum_i\op A_i\op H^i$. Hence
$$
 J=\{\op A\in\mathcal D:[\op A]\in W\}
$$
is a nonzero proper left ideal, stable under right multiplication by
$\op H$ because $[\op A\op H]=\eta[\op A]$. Write $J=\mathcal D\op A$.
Since $\op A\op H\in\mathcal D\op A$, $\op H$ preserves the local kernel of $\op A$.
Section~\ref{sec:preliminaries} gives $f(\op H)\in J$ for some
$0\ne f\in\mathbb C[z]$. But
$[f(\op H)]=f(\eta)[1]$, with $f(\eta)$ invertible, so $[1]\in W$,
a contradiction. Differential division gives dimension $k$.
\end{proof}

\begin{theorem}
\label{thm:order-product-identity}
Let a polynomial differential operator of order $r$ have monic
eigenpolynomials in every sufficiently large degree, with eigenvalue
polynomial $\lambda$ of degree $d\ge1$ and a fixed recurrence of
order $m$. If $m\ne2$, then $r=dm$. The same equality holds for
$m=2$ if no $\zeta\in\mathbb C$ makes $\lambda(\zeta+z)$ even.
\end{theorem}

\begin{proof}
For $m=1$, the proof of Theorem~\ref{thm:spectral-degree-divisibility} gives $r=d$.
For $m\ge2$, its proof gives $\op N\in\mathcal Q$ with
$\lambda(\op N)=\op L$ and diagonal class $\op N_{\rm d}P_n(x)=nP_n(x)$.
Take $\op T,\op H,\op R,F$ from
Lemma~\ref{lem:ordinary-root-normalization}, with $k=\ord\op H=r/d$.
Cancellation also gives
$$
 \op N\op R=\op R\op H.
$$

The unweighted formal adjoint
$\op A^*=\sum_i(-\partial)^i a_i(x)$ of
$\op A=\sum_i a_i(x)\partial^i$ is $\mathbb C$-linear and reverses
products. It differs from transposition on functionals.

Put $K=\mathbb C(\eta)$, where $\eta$ is an indeterminate, and let
$U=\bigoplus_{j\in\mathbb Z}K e_j$.  With $w=m-1$, define
$$
 x e_j=e_{j+1}+\sum_{s=0}^w b_s(\eta+j)e_{j-s},\quad
 \op L e_j=\lambda(\eta+j)e_j.
$$
Let $\mathcal B=\mathbb C\langle x,\op L\rangle$ be the unital
algebra generated by $x$ and $\op L$. Its identities hold on $U$:
their finite rational difference representations vanish on the original
polynomial tail, hence identically, before substituting $n=\eta+j$.
Let
$$
 M=\bigoplus_{j\in\mathbb Z}K\mathbf a_j\subseteq U^*,\quad
 \langle\mathbf a_j,e_i\rangle=\delta_{ji},
$$
be its restricted dual, where $U^*=\operatorname{Hom}_K(U,K)$.
The spaces $U$ and $M$ are algebraic direct sums over $K$;
the pairing precedes localization in $x$. Every matrix from $\mathcal B$
has finitely many nonzero diagonals, so its transpose preserves $M$:

$$
 \begin{aligned}
 x'\mathbf a_j&=\mathbf a_{j-1}
   +\sum_{s=0}^w b_s(\eta+j+s)\mathbf a_{j+s},\\[7pt]
 \op L'\mathbf a_j&=\lambda(\eta+j)\mathbf a_j.
 \end{aligned}
$$
Primes refer to this pairing. Letting $\op A^*$ act as $\op A'$
defines a left action of $\mathcal B^*=\mathbb C\langle x,\op L^*\rangle$:
the two product reversals cancel. This is the algebraic convention
$\Dd=-\partial'$ of Section~\ref{sec:preliminaries}.
For $\op C\in K\otimes_{\mathbb C}\mathcal B^*$ of order $q$ and
$0\ne f\in K[x]$, commutation with $f$ lowers differential order,
so $(\ad_f)^{q+1}\op C=0$. The binomial identities give
$$
 \begin{aligned}
 f^{q+1}\op C
 &=\sum_{i=0}^{q}\binom{q+1}{i}
   (\ad_f)^i(\op C)f^{q+1-i},\\[7pt]
 \op C f^{q+1}
 &=\sum_{i=0}^{q}(-1)^i\binom{q+1}{i}
   f^{q+1-i}(\ad_f)^i(\op C).
 \end{aligned}
$$
These are a right and a left multiple of $f$, respectively, within
$K\otimes_{\mathbb C}\mathcal B^*$. They allow localization by
the nonzero polynomials in $K[x]$.
Also $(\ad_x)^{r-1}\op L^*=a(x)\partial+b(x)$ with $a\ne0$;
coefficient inversion recovers
$\partial=a^{-1}((\ad_x)^{r-1}\op L^*-b)$. Thus the coefficient
localization is $\mathcal D_\eta$. Over $\mathbb C$, the same argument
localizes $\mathcal B$ to $\mathcal D$, with common right coefficient
denominators.

The recurrence generates $M$ over $K[x]$ from
$\mathbf a_0,\ldots,\mathbf a_w$ by solving in both directions.
For $q\ge0$, $x^q\mathbf a_i$ has extreme indices $i-q,i+wq$,
with nonzero coefficients. In a nontrivial relation
$\sum_{i=0}^w f_i(x)\mathbf a_i=0$, put $D=\max_i\deg f_i$
and choose the largest $i$ with $\deg f_i=D$.
If $i>0$, the largest index $i+wD$ is unique.
If $i=0$, the smallest index $-D$ is unique.
Neither coefficient can vanish. Thus $M$ is free on these generators,
and coefficient localization is injective. Hence
$$
 W=K(x)\otimes_{K[x]}M
$$
is a left $\mathcal D_\eta$-module of dimension $m$ over $K(x)$,
generated by $\mathbf a_0$. Indeed, the values $\lambda(\eta+j)$ are
distinct because a nonconstant polynomial has no nonzero additive period.
Interpolation in $\op L^*$ on each finite transposed recurrence isolates
the nonzero extreme steps $-1,+w$. Starting at zero, these steps reach every integer,
hence every coordinate functional.

Choose $0\ne p\in\mathbb C[x]$ with $\op R p,\op R\op H p\in\mathcal B$.
Since $\op N\op R=\op R\op H$, the finite rational bands of
$\op R\op H p$ and $\op N_{\rm d}(\op R p)$ differ by finite rank,
hence coincide. Thus on $U$,
$$
 (\op R\op H p)=\op N_\eta(\op R p),\quad
 \op N_\eta e_j=(\eta+j)e_j.
$$
Apply the transposed identity to $\mathbf a_0$ and pass to $W$:
$$
 p\op H^*\op R^*\mathbf a_0
   =p\op R^*\op N_\eta'\mathbf a_0
   =\eta p\op R^*\mathbf a_0.
$$
Cancel $p(x)$ to obtain
$$
 \op H^*q=\eta q,\quad q=\op R^*\mathbf a_0.
$$
Furthermore,
$$
 \op T^*q=F(\op L^*)\mathbf a_0
    =F(\lambda(\eta))\mathbf a_0\ne0.
$$
Thus $q\ne0$ and $\mathbf a_0\in\mathcal D_\eta q$, so
$W=\mathcal D_\eta q$. The map $[\op A]\mapsto\op A q$ is a nonzero
surjection
$$
 \mathcal D_\eta/\mathcal D_\eta(\op H^*-\eta)\longrightarrow W.
$$
Its source is simple of dimension $k$ by
Lemma~\ref{lem:generic-eigenvalue-simple}. Hence it is an isomorphism
and $m=k=r/d$.
\end{proof}

For a nonzero rational function $f(x)=p(x)/q(x)$, write
$\deg_\infty f=\deg p-\deg q$. At a finite point $a$,
$\operatorname{ord}_a f$ is its order of vanishing; a pole has negative
order. Set $\deg_\infty0=-\infty$ and $\operatorname{ord}_a0=+\infty$.

\begin{proposition}
\label{prop:general-bispectral-reduction}
Let $\op L$ be a nonscalar polynomial differential eigenoperator for
$(P_n(x))_{n\ge N}$ with a recurrence of fixed length.
Suppose that $h,\mu\in\mathbb C[z]$ are nonconstant and
$h(\op N)=\op L$ for some $\op N\in\mathcal Q$ whose class
modulo finite rank is diagonal with entries $\mu(n)$.
There are a rational differential operator $\op H_0$ of order $\ord\op N$,
rational functions $\phi_n(x)$, and a nonzero rational differential
operator $\op Q$ such that
$$
 \begin{gathered}
 \op H_0\phi_n(x)=\mu(n)\phi_n(x),\quad
 \op Q\phi_n(x)=f(n)P_n(x),\\[7pt]
 \op L\op Q=\op Qh(\op H_0),
 \end{gathered}
$$
where $f(z)$ is a nonconstant polynomial. The family
$(\phi_n(x))_{n\ge N}$ is bispectral, and a finite rational bispectral Darboux transformation
recovers the prescribed polynomials after division by $f(n)$.
All identities hold in every sufficiently large degree.
\end{proposition}

\begin{proof}
Apply Lemma~\ref{lem:ordinary-root-normalization} to $h(\op N)=\op L$,
obtaining $\op T,\op R,\op H,F$. Here $\lambda=h\circ\mu$.
The coefficient localization of $\mathcal B=\mathbb C\langle x,\op L\rangle$
is $\mathcal D$ on both sides, as proved above. Choose polynomials
$p,q\ne0$ with
$$
 \op P=p(x)\op T\in\mathcal B,\quad
 \op Q=\op R q(x)\in\mathcal B,
$$
and $p\op H\op T\in\mathcal B$. Set $\theta=pq$ and
$f=F\circ\lambda$, taking $p\mapsto xp$ if needed to make $\theta$
nonconstant. Put
$$
 \phi_n(x)=q(x)^{-1}\op T P_n(x),\quad
 \op H_0=q(x)^{-1}\op H q(x).
$$
For $\op M_\mu P_n(x)=\mu(n)P_n(x)$, the map
$p\op H\op T-(p\op T)\op M_\mu$ has finite rank and a finite rational
difference representation, so vanishes on the tail by the preceding argument. Thus
$\op H\op T P_n(x)=\mu(n)\op T P_n(x)$, giving the source equation.
The differential identities give
$$
 \begin{gathered}
 F(\op L)=\op Q\theta(x)^{-1}\op P,\quad
 \phi(x)=\theta(x)^{-1}\op P\psi(x),\\[7pt]
 \op Q\phi(x)=f(n)\psi(x),\quad
 \op L\op Q=\op Qh(\op H_0).
 \end{gathered}
$$
Here $\psi(x)=(P_n(x))_{n\ge N}$ and
$\phi(x)=(\phi_n(x))_{n\ge N}$ are columns. Increase $N$ so that
$f(n)\ne0$ for $n\ge N$, and let $(\mathbf a_j)_{j\ge0}$ be the dual of a monic completion
of $(P_n(x))_{n\ge N}$. The functionals on the span of $(\phi_n(x))_{n\ge N}$ given by
$$
 \ell_j(v)=f(j)^{-1}\duality{\mathbf a_j}{\op Qv},\quad
 \ell_j(\phi_n(x))=\delta_{jn},\quad j,n\ge N,
$$
prove the independence of $(\phi_n(x))_{n\ge N}$ and give the coordinates for $\op b_\phi^{-1}$.

For the reverse step put $\op U=\op b^{-1}(\op P)$ and $\op V=\op b^{-1}(\op Q)$,
using the original family's inverse Fourier map. Both are finite rational difference operators.
Reorder the factorization in $\mathcal Q$ as
$\theta(x)=\op P F(\op L)^{-1}\op Q$.
Modulo finite rank, the inverse is diagonal with entries $f(n)^{-1}$.
The coordinate criterion gives
$$
 \op b^{-1}(\theta)=\op V f(n)^{-1}\op U.
$$
Put $\op S=\theta^{-1}\op P\theta$. Since
$\theta(x)\phi(x)=\op U\psi(x)$, the preceding identity gives
$\psi(x)=(\op V f(n)^{-1})\phi(x)$. Using the coordinate functionals
$\ell_j$ gives
$$
 \begin{aligned}
 \op b_\phi^{-1}(\op Q)&=f(n)\op V f(n)^{-1},\\[7pt]
 \op b_\phi^{-1}(\op S)&=\op U,\\[7pt]
 \op b_\phi^{-1}(\theta)&=\op U\op V f(n)^{-1}.
 \end{aligned}
$$
Thus the reverse factors and their multiplier lie in the same Fourier
algebra, including when $\phi_n(x)$ has poles.
The identities
$$
 \begin{gathered}
 F(h(\op H_0))
 =\theta(x)^{-1}\op P\op Q
 =\op S\theta(x)^{-1}\op Q,\\[7pt]
 \op Q\op S\theta(x)^{-1}=F(\op L)
 \end{gathered}
$$
give an admissible reverse bispectral Darboux step to $f(n)P_n(x)$.
Monic normalization removes $f(n)$, nonzero for all large $n$.
The identity $\op L\op Q=\op Qh(\op H_0)$ identifies $\op L$
in the eigenoperator algebra of the transformed family.

The eigenvalues $\mu(n),\theta(x)$ are nonconstant.
Zeros and poles of rational coefficients in $n$ exclude finitely many indices.
\end{proof}

The rational source has a difference equation for $\theta(x)$.
For recurrence order at least three,
Theorem~\ref{thm:all-polynomial-sources} supplies a polynomial source
with a recurrence for $x$.

The bispectrality of the two Hermite Darboux transforms in
\cite[Sections~4--5, pp.~928--935]{castrogrunbaum2026} follows, on a tail,
from the same argument with the dual basis.
The polynomial recurrences determine finite-band identities,
which act on both independent source solutions by the Hermite recurrence
and its derivative \cite[Section~3, pp.~925--927]{castrogrunbaum2026}.
The target operators have rational coefficients and polynomial
multipliers $x^2/2$ and $2x^3/3-x$.

\section{Affine eigenvalues}
\label{sec:affine-eigenvalues}

\subsection{The differential algebra}

\begin{proposition}
\label{prop:affine-spectral-algebra}
Let $\op L$ have polynomial coefficients and monic eigenpolynomials
$(P_n(x))_{n\ge0}$ with eigenvalues $n$.
The differential operators with rational coefficients commuting with
$\op L$ form $\mathbb C[\op L]$. For every $N\ge0$, this is also the
algebra of rational differential operators having $(P_n(x))_{n\ge N}$
as eigenpolynomials. Their minimum nonscalar order is $\ord\op L$.
In the expansion
$$
 xP_n(x)=P_{n+1}(x)+\sum_{s=0}^n b_s(n)P_{n-s}(x),
$$
each $b_s(z)$ is polynomial and divisible by $z^{\underline s}$.
Consequently, a recurrence of fixed length valid in sufficiently large
degrees holds for this eigenpolynomial sequence in every degree.
\end{proposition}

\begin{proof}
Write
$$
 \op L=\op E+\sum_{i=1}^h x^{-i}p_i(\op E),
 \quad z^{\underline i}\mid p_i(z).
$$
For $n\ge0$, if $0\ne f\in\mathbb C(x)$ satisfies $\op Lf=nf$, its leading term
$cx^m$ at infinity gives $m=n$. Applying this to $f-cP_n(x)$ yields
$f=cP_n(x)$.

If $\op A\in\mathcal D$ commutes with $\op L$, then
$\op AP_n(x)=\mu(n)P_n(x)$. Thus $\op A\Pp\subset\Pp$, so $\op A$ has
polynomial coefficients. Comparison of leading terms shows that $\mu$ is polynomial, and
$\op A-\mu(\op L)$ annihilates a basis. An operator diagonal only
in high degrees also commutes with $\op L$: the commutator has
infinite-dimensional kernel and is zero. For nonconstant $\mu$,
$\ord\mu(\op L)=(\deg\mu)\ord\op L$.

Put
$P_n(x)=\sum_{j=0}^n c_j(n)x^{n-j}$. Comparison of coefficients gives
$$
 \begin{aligned}
 jc_j(z)&=\sum_{i=1}^{\min(j,h)}p_i(z-j+i)c_{j-i}(z),\quad c_0(z)=1,\\[7pt]
 b_s(z)&=c_{s+1}(z)-c_{s+1}(z+1)
 -\sum_{t=0}^{s-1}b_t(z)c_{s-t}(z-t).
 \end{aligned}
$$
These formulas show by induction that $c_j(z)$ and $b_s(z)$ are
polynomials, with $z^{\underline j}\mid c_j(z)$ and
$z^{\underline s}\mid b_s(z)$.

By Lemma~\ref{lem:dual-spectral-support}, the eventual recurrence
holds on the full dual sequence and hence in every degree.
\end{proof}

\subsection{Commuting components}

The homogeneous component of degree change $t$ is the sum of the
terms $x^{j+t}\partial^j$; it has the form $x^t p(\op E)$.
For rational coefficients, use their Laurent expansions at infinity.

For operators commuting with a fixed homogeneous operator, see
\cite{bavula2005}.
For polynomial generators, Horozov gives the construction, raising and
lowering operators, and recurrence \cite[Theorems~4.1--4.2]{horozov2018vector}.

\begin{proposition}
\label{prop:commuting-recurrence}
Let $S$ be a finite nonempty set of positive integers and let
$$
 \op{A}_i=x^{-i}p_i(\op{E}),\quad i\in S,
$$
be nonzero differential operators with polynomial coefficients,
where $p_i(z)$ are polynomials. Suppose that the $\op{A}_i$ commute
pairwise, and put
$$
 \op{L}=\op{E}+\sum_{i\in S}\op{A}_i,\quad
 \op{Q}=\sum_{i\in S}\frac{\op{A}_i}{i}.
$$
The unique monic eigenpolynomials are $P_n(x)=e^{\op{Q}}x^n$,
with $\op{L}P_n(x)=nP_n(x)$ for $n\ge0$. They satisfy a recurrence
of fixed length in sufficiently large degrees exactly when there
are a nonconstant polynomial $g(z)$ and nonzero constants $c_i$ such that
\refstepcounter{equation}\label{eq:commuting-root-criterion}
$$
 \op{G}=x^{-1}g(\op{E}),\quad \op{A}_i=c_i\op{G}^i,\quad i\in S.
 \leqno\textup{(\theequation)}
$$
In this case, with $K=\max S$ and $\sigma=\deg g$, the lowest nonzero recurrence
shift is $1-K\sigma$. For every $N\ge0$, the differential
operators with rational coefficients having $(P_n(x))_{n\ge N}$
as eigenpolynomials form $\mathbb C[\op{L}]$, and their minimum
nonscalar order is $K\sigma$.
\end{proposition}

\begin{proof}
Since
$[\op{Q},\op{E}]=\sum_i\op{A}_i$, conjugation gives
$\op{L}=e^{\op{Q}}\op{E}e^{-\op{Q}}$. Each exponential is finite
on every polynomial. For the monomial dual
$\mathbf e_n=(-1)^n\Dd^n\boldsymbol\delta_0/n!$, the dual sequence is
$\mathbf a_n=(e^{-\op Q})'\mathbf e_n$.
The exponential series converges in $\sigma(\Pp',\Pp)$, since
$\op Q^jp=0$ for large $j$ for every $p\in\Pp$.
Write
$$
 \op J=e^{-\op Q}xe^{\op Q},\quad
 \op V=e^{-\op Q}\op E e^{\op Q}=\op E-\sum_i\op A_i.
$$
Then
$$
 \duality{\mathbf e_m}{\op Jx^n}
 =\duality{\mathbf a_m}{xP_n(x)}.
$$
An eventual recurrence holds globally by
Lemma~\ref{lem:dual-spectral-support}; uniqueness of weak dual
coordinates therefore makes $\op J$ a finite homogeneous sum.

Consider formal sums
$\sum_{t\le M}x^t f_t(\op E)$, where $f_t\in\mathbb C(z)$ and
$f(\op E)x=xf(\op E+1)$. Use the degree filtration: a remainder
tends to zero when its largest degree tends to $-\infty$.
On sums bounded above in degree, every fixed degree in a product
receives only finitely many terms. Since $\op Q$ strictly lowers
degree, the same holds for its conjugation series.
In these symbolic calculations every nonzero $x^t f(\op E)$ has
an inverse; the rational functions of $\op E$ are kept as symbols.
Here $e^{\ad_{\op Q}}\op A=e^{\op Q}\op A e^{-\op Q}$.
Put $\delta=\gcd S$ and choose
integers $u_i$ with $\sum_i iu_i=\delta$. The commuting product
$\op{H}=\prod_i\op{A}_i^{u_i}$ has degree $-\delta$.
Each $\op A_i\op H^{-i/\delta}$ has degree zero, so it is a
rational function of $\op E$. Since it commutes with $\op H$,
that rational function has period $\delta$ and must be constant.
Thus $\op A_i=d_i\op H^{i/\delta}$.
Let $\op{J}_N$ be the lowest nonzero component of $\op{J}$, of
degree $1-\delta N$. The lowest component of
$[\op{V},\op{J}]=\op{J}$ gives
$[\op{H}^{K/\delta},\op{J}_N]=0$, and hence
$[\op{H},\op{J}_N]=0$. Indeed, if $\op{C}=[\op{H},\op{J}_N]$
were nonzero, all terms in
$\sum_{r=0}^{K/\delta-1}\op{H}^r\op{C}\op{H}^{K/\delta-1-r}$
would have the same nonzero leading term in their rational
coefficients in $\op E$ at infinity.
It follows that
$\op G=\op H^N\op J_N^{-1}=x^{-1}g(\op E)$ commutes with all
$\op A_i$. The degree-zero ratio $\op A_i\op G^{-i}$ is then a
rational function of $\op E$ with period one, and is constant.
This gives $\op A_i=c_i\op G^i$ with $g$ rational.

Write $g=a/b$ in lowest terms and
$Q(z)=\sum_i c_i z^i/i$.
Conjugating $x\op{G}b(\op{E})=a(\op{E})$ degree by degree gives
$$
 \op{J}\op{G}b(\op{V})=a(\op{V}),\quad
 \op{V}=\op{E}-\op{G}Q'(\op{G}).
$$
All expressions are finite sums of powers of $\op{G}$ with
rational coefficients in $\op{E}$. Choose polynomials $r(z),s(z)$
with $r(z)a(z)+s(z)b(z)=1$. Then
$$
 (r(\op{V})\op{J}\op{G}+s(\op{V}))b(\op{V})=\op{I}.
$$
In a product, the largest exponents add, as do the smallest.
A product equal to the identity therefore forces $b(\op V)$ to
contain only one power of $\op G$. If $\deg b>0$, however,
$b(\op V)$ has nonzero coefficients at powers $0$ and $K\deg b$.
Thus $b$ is constant and $g$ is polynomial. The operators
$\op A_i=c_i\op G^i$ could not have polynomial coefficients if
$g$ were constant, so $g$ is nonconstant.

Conversely, put $\Delta_i f(z)=f(z+i)-f(z)$. For $f\in\mathbb C[z]$ and $m\in\mathbb Z$, the symbolic algebra gives
$$
 [\op G^i,f(\op E)\op G^m]
 =\Delta_i f(\op E)\op G^{m+i}.
$$
Since $x=g(\op E)\op G^{-1}$, each commutator lowers the polynomial
degree in $\op E$; after the first, all exponents of $\op G$ are
nonnegative. Hence $(\ad_{\op Q})^{\sigma+1}x=0$, and $\op J$ has
degrees between $1-K\sigma$ and $1$. Here
$\op Q=\sum_i\op A_i/i$ has polynomial coefficients and lowers degree,
so this symbolic conjugation equals its locally finite action on $\Pp$.
The dual coordinates above give the full recurrence.
The lowest component is
$$
 \frac{(-c_K/K)^\sigma}{\sigma!}
 \Delta_K^\sigma g(\op E)\op G^{K\sigma-1}
 =(-c_K)^\sigma[z^\sigma]g(z)\op G^{K\sigma-1}.
$$
Here $\Delta_K^\sigma g=\sigma!K^\sigma[z^\sigma]g$.
Acting on $x^n$ gives
\refstepcounter{equation}\label{eq:commuting-last-coefficient}
$$
 b_{K\sigma-1}(n)=(-c_K)^\sigma[z^\sigma]g(z)
 \prod_{j=0}^{K\sigma-2}g(n-j).
 \leqno\textup{(\theequation)}
$$
This polynomial is nonzero.
The algebra and minimum-order assertions follow from
Proposition~\ref{prop:affine-spectral-algebra}, since
$\ord\op L=K\sigma$.
\end{proof}

The identity
$$
 \op{G}^i=x^{-i}\prod_{r=0}^{i-1}g(\op{E}-r)
$$
shows that $\op{G}^i$ has polynomial coefficients exactly when
$z^{\underline i}$ divides $\prod_{r=0}^{i-1}g(z-r)$.
Write $Q(z)=\sum_{i\in S}c_i z^i/i$, so $\op Q=Q(\op G)$.
For $g(0)\ne0$, a positive integer root $b$ also gives a
representation with a polynomial generator at the original degree indices.
Put
$$
 \widetilde g(z)=\frac{zg(z)}{z-b},\quad
 \widetilde{\op G}=x^{-1}\widetilde g(\op E),\quad
 \op F=\prod_{j=0}^{b-1}(\op E-j)=x^b\partial^b.
$$
Then $\widetilde{\op G}$ has polynomial coefficients and
$\op F\widetilde{\op G}=\op G\op F$. Consequently
$$
 P_n(x)=\frac{\op F e^{Q(\widetilde{\op G})}x^n}{n^{\underline b}},
 \quad n\ge b.
$$
When $g(0)=0$, the operator
$\op{L}=\op{E}+Q'(\op{G})\op{G}$ has the second form listed in
\cite[Conjecture~1.10]{hst2025}; the first form is the case $S=\{1\}$.
Theorem~\ref{thm:affine-global-normal-form} below proves that these
components commute after a translation whenever the eigenvalues
are affine and a recurrence of fixed length is present.

\section{Classification of irreducible operators}
\label{sec:principal-low-order}

Give $x,\partial$ rational weights $u,v$, with $u+v>0$;
$x^a\partial^j$ has weight $ua+vj$.
Write $\mathrm{wt}(\op A)$ for the largest weight of its terms.
The weighted symbol is the sum of those terms, with $\partial$
replaced by a commuting variable $\xi$. Their exponent pairs lie on
a supporting face of the convex hull of all exponent pairs present.
For weights $(0,1)$ this is the principal symbol
$\operatorname{lc}(\op A)\xi^{\ord\op A}$.
For rational coefficients and $u>0$, use their Laurent expansions at infinity.

Following \cite[Lemma~2.7]{dixmier1968}, define the derivation
$\op H_h$ and the Poisson bracket by
$$
 \op H_h=h_\xi\partial_x-h_x\partial_\xi,\quad
 \{h,f\}=\op H_h f.
$$
Along $\dot x=h_\xi$, $\dot\xi=-h_x$, one has $\dot f=\op H_hf$.

\begin{lemma}
\label{lem:nonlinear-weighted-face}
Suppose $\op L$ has monic eigenpolynomials with polynomial eigenvalue
$\lambda(n)$ of degree $d\ge1$ and a recurrence of fixed length in
sufficiently large degrees. Give $x,\partial$ rational weights $u,v$
with $\delta=u+v>0$. If the weighted leading expression $h$ of
$\op L$ has weight $\eta>d\delta$, then
$\op H_h^N x=0$ for some $N\ge1$.
\end{lemma}

\begin{proof}
Let $J$ be the set of shifts with nonzero recurrence coefficients. Use the commuting maps $\Delta,\mathcal S$ of
Lemma~\ref{lem:dual-spectral-support}, and put
$h^\vee(x,\zeta)=h(x,-\zeta)$.
Since $\partial'=-\Dd$ and reordering lowers weight by at least $\delta$,
transposition preserves weights, with leading expression $h^\vee$.
Thus $\Delta$ raises weight by at most $\eta-\delta$, with
leading action $-\op H_{h^\vee}$, whereas $\mathcal S$ raises it
by at most $\eta$. For $1\le k\le d$, the term
$c_{jk}(\mathcal S)\Delta^{d-k}$ raises weight by at most
$$
 (d-1)\eta-(d-k)\delta<d(\eta-\delta).
$$
Hence, with $N=d|J|$, a nonzero $(-\op H_{h^\vee})^Nx$ would
be the leading expression of the zero dual operator
$\prod_{j\in J}R_j(\Delta,\mathcal S)(x)$, a contradiction.
The substitution $\zeta=-\xi$ identifies it with $\op H_h^Nx$.
\end{proof}

\begin{lemma}
\label{lem:formal-face-obstruction}
Let $z=z(x,\xi)$ be a local coordinate, with rational powers of $x$ allowed. Suppose that
$$
 h=x^d z^k H(z),\quad z_\xi=x^b B(z),\quad
 d,b\in\mathbb Q,\quad k\in\mathbb Z_{>0},\quad H(0)B(0)=C\ne0,
$$
where $H(z),B(z)\in\mathbb C[[z]]$.
If $e=b+d/k$ satisfies $e<0$ or $e>1$, then
$\op H_h^N x\ne0$ for every $N\ge1$.
\end{lemma}

\begin{proof}
In these coordinates,
$$
 \op H_hx=Ck x^{d+b}z^{k-1}+O(z^k),\quad
 \op H_hz=-Cd x^{d+b-1}z^k+O(z^{k+1}).
$$
Induction gives
$$
 \op H_h^N x=C^N
 \prod_{j=0}^{N-1}\bigl(k+j(d+k(b-1))\bigr)
 x^{1+N(d+b-1)}z^{N(k-1)}+O(z^{N(k-1)+1}).
$$
Each factor is $k(1+j(e-1))$. If $e>1$, all are positive.
If $e<0$, the only possible zero has
$j=1/(1-e)\in(0,1)$, so no factor vanishes.
\end{proof}

\begin{theorem}
\label{thm:general-newton-triangle}
Suppose $\op L$ has polynomial coefficients, differential order $r$,
polynomial eigenvalue $\lambda(n)$ of degree $1\le d<r$, and monic
eigenpolynomials with a recurrence of fixed length in sufficiently
large degrees. After translating $x$, its principal coefficient and
homogeneous decomposition have the following form, where $\sigma$
is a positive divisor of $r$:
$$
 \begin{aligned}
 a_r(x)&=cx^{r-r/\sigma},\quad c\ne0,\quad\sigma\mid r,
 \quad K=r/\sigma,\\[7pt]
 \op L&=\lambda(\op E)+\sum_{i=1}^K x^{-i}p_i(\op E),\quad
 \deg p_i\le\left\lfloor d+\frac{r-d}{K}i\right\rfloor,
 \quad\deg p_K=r.
 \end{aligned}
$$
No commutation hypothesis is imposed on the degree-decreasing terms.
\end{theorem}

\begin{proof}
Use weights $(0,1)$.
Since $r>d$, Lemma~\ref{lem:nonlinear-weighted-face} gives
$\op H_h^N x=0$ for some $N$, where $h=a_r(x)\xi^r$.
Along a local solution of the associated differential equations,
the Taylor series of the first coordinate $X(t)$ therefore
terminates. Thus $X(t)$ is polynomial. On $h=1$ it satisfies
$$
 (X')^r=r^r a_r(X).
$$
Since $d<r$, one has $\deg a_r<r$. If a root of $a_r$ has
multiplicity $m$ and a preimage under the nonconstant polynomial $X$
has multiplicity $q$, then $r(q-1)=mq$, so $q\ge2$ and $m\ge r/2$.
There can therefore be only one distinct root. Its multiplicity is
$r-r/\sigma$ with $\sigma\mid r$; the constant case is $\sigma=1$.
Translate the root when present and put $K=r/\sigma$.

Write
$\op L=\lambda(\op E)+\sum_{i\ge1}x^{-i}p_i(\op E)$.
Represent a term $x^a\partial^j$ by the pair $(i,j)$, where
$i=j-a$ is its decrease in degree. The principal term corresponds
to $(K,r)$.

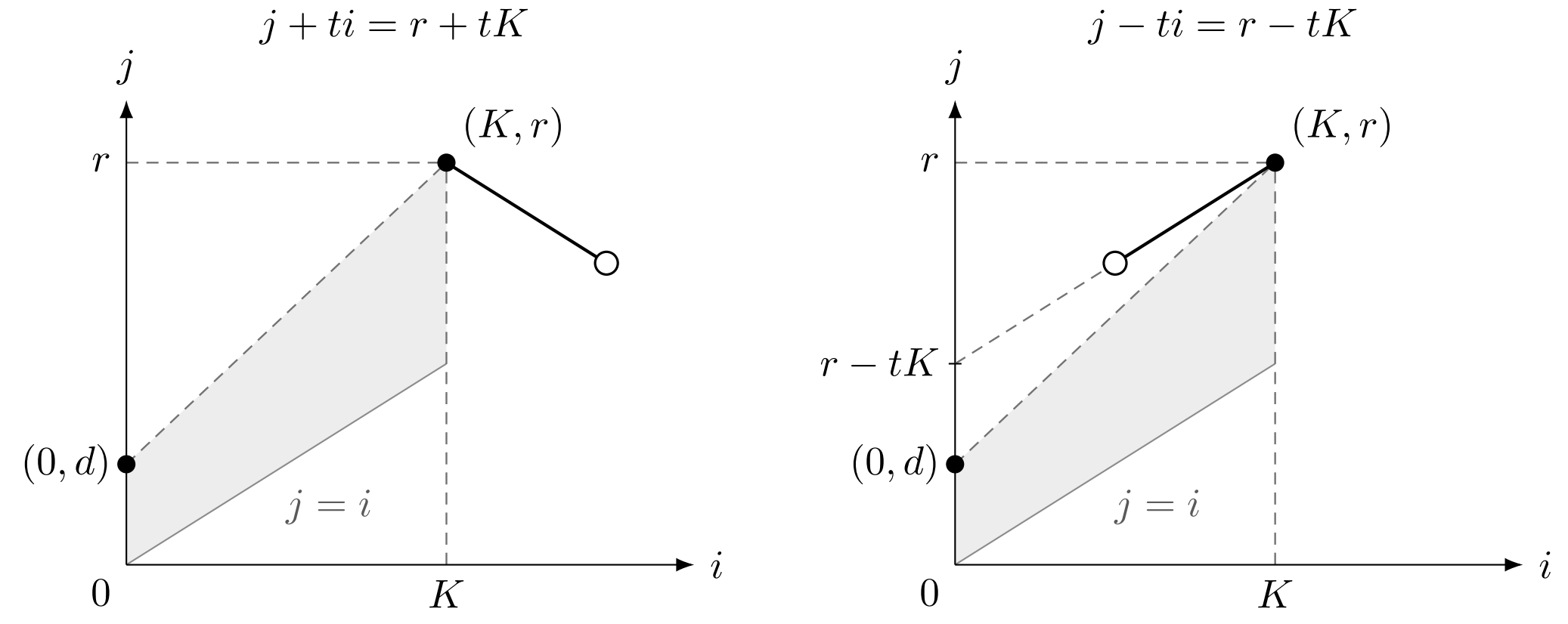
\begin{figure}[ht]
\centering
\begin{tikzpicture}[x=.70cm,y=.44cm,font=\small,
  axis/.style={-{Latex[length=1.6mm]},line width=.4pt},
  bound/.style={densely dashed,black!55,line width=.45pt},
  face/.style={line width=.8pt},
  support point/.style={circle,fill=black,inner sep=1.6pt},
  assumed point/.style={circle,draw=black,fill=white,line width=.6pt,
    inner sep=2pt}]
\begin{scope}
  \fill[black!7] (0,0)--(0,2)--(4,8)--(4,4)--cycle;
  \draw[black!45,line width=.4pt] (0,0)--(4,4);
  \draw[bound] (0,2)--(4,8);
  \draw[bound] (4,0)--(4,8);
  \draw[bound] (0,8)--(4,8);
  \draw[axis] (0,0)--(7.1,0) node[right] {$i$};
  \draw[axis] (0,0)--(0,9.25) node[above] {$j$};
  \node[below left] at (0,0) {$0$};
  \node[below] at (4,0) {$K$};
  \node[left] at (0,8) {$r$};
  \node[support point] at (0,2) {};
  \node[left] at (0,2) {$(0,d)$};
  \node[support point] at (4,8) {};
  \node[above right] at (4,8) {$(K,r)$};
  \node[below right,black!65] at (1.8,1.8) {$j=i$};
  \draw[face] (4,8)--(6,6);
  \node[assumed point] at (6,6) {};
  \node[above] at (3.35,10.05) {$j+ti=r+tK$};
  \node at (3.35,-2.0) {\textup{(a)} $i>K$};
\end{scope}
\begin{scope}[xshift=7.25cm]
  \fill[black!7] (0,0)--(0,2)--(4,8)--(4,4)--cycle;
  \draw[black!45,line width=.4pt] (0,0)--(4,4);
  \draw[bound] (0,2)--(4,8);
  \draw[bound] (4,0)--(4,8);
  \draw[bound] (0,8)--(4,8);
  \draw[axis] (0,0)--(7.1,0) node[right] {$i$};
  \draw[axis] (0,0)--(0,9.25) node[above] {$j$};
  \node[below left] at (0,0) {$0$};
  \node[below] at (4,0) {$K$};
  \node[left] at (0,8) {$r$};
  \node[support point] at (0,2) {};
  \node[left] at (0,2) {$(0,d)$};
  \node[support point] at (4,8) {};
  \node[above right] at (4,8) {$(K,r)$};
  \node[below right,black!65] at (1.8,1.8) {$j=i$};
  \draw[bound] (0,4)--(2,6);
  \draw[face] (2,6)--(4,8);
  \node[assumed point] at (2,6) {};
  \draw[line width=.4pt] (-.08,4)--(.08,4);
  \node[left] at (0,4) {$r-tK$};
  \node[above] at (3.35,10.05) {$j-ti=r-tK$};
  \node at (3.35,-2.0) {\textup{(b)} $j>d+(r-d)i/K$};
\end{scope}
\end{tikzpicture}
\caption{Supporting faces in the proof of
Theorem~\ref{thm:general-newton-triangle}, drawn schematically for $1<K<r$.
The shaded region is $0\le i\le K$, $i\le j\le d+(r-d)i/K$.
Open circles represent the hypothetical terms excluded in the two cases.
The supporting parameter satisfies $t>0$ in \textup{(a)} and
$0<t<(r-d)/K$ in \textup{(b)}.}
\label{fig:newton-support-faces}
\end{figure}

Suppose first that a term has $i>K$, as in
Figure~\ref{fig:newton-support-faces}\textup{(a)}. For weights $(-t,1+t)$,
increase $t$ from zero until the principal term first has the same
weight as another term. This first value is positive, and the sum
of terms of that weight has the form
$$
 \begin{aligned}
 h&=x^\alpha\xi^\beta F(x^p\xi^q),\quad
 p,q\in\mathbb Z_{>0},\quad p>q,\quad F(0)\ne0,\quad M=\deg F\ge1,\\[7pt]
 \alpha+pM&=r-K,\quad\beta+qM=r,\quad p/q=(1+t)/t.
 \end{aligned}
$$
Its weight is $r+tK>d$. At a nonzero root $z_0$ of multiplicity
$k\le M$, put $z=x^p\xi^q-z_0$. Then
$$
 h=x^{d_0}z^k H(z),\quad z_\xi=x^{p/q}B(z),\quad
 d_0=-K-r/t,\quad H(0)B(0)\ne0.
$$
The exponent in Lemma~\ref{lem:formal-face-obstruction} is
$$
 e=1-\frac{K+(r-k)/t}{k}<0,
$$
because $kp\le Mp\le r-K<r(p-q)+Kq$.
This contradicts Lemma~\ref{lem:nonlinear-weighted-face}.
Thus no degree decrease exceeds $K$.

If a remaining monomial satisfies $j>d+(r-d)i/K$, then $i<K$;
see Figure~\ref{fig:newton-support-faces}\textup{(b)}.
For weights $(t,1-t)$ take the first supporting face through $(K,r)$:
$$
 0<t=\min_{1\le i<K}\frac{r-j}{K-i}<\frac{r-d}{K},
$$
where the minimum is taken over the terms present in $\op L$ with
$1\le i<K$. Its weight $r-tK$ is greater than $d$, so the entire diagonal
part $\lambda(\op E)$ lies strictly below it.
Write $t=p/q$ in lowest terms, with $p,q>0$.
The face indices are $(K-qm,r-pm)$, $0\le m\le M$.
Since they all have positive degree decrease, $qM\le K-1$, and
$$
 h=x^{r-K}\xi^r F(Z),\quad Z=x^{q-p}\xi^{-p},\quad
 F(0)\ne0,\quad\deg F=M\ge1.
$$
At a nonzero root $z_0$ of multiplicity $k\le M$, put $z=Z-z_0$.
Then
$$
 h=x^{d_0}z^kH(z),\quad z_\xi=x^{1-1/t}B(z),\quad
 d_0=r/t-K>0,\quad H(0)B(0)\ne0.
$$
Here $e=1-1/t+d_0/k>1$, since
$$
 qr-Kp=K(\sigma q-p)\ge K>qM\ge qk.
$$
The integer $\sigma q-p$ is positive because $t<(r-d)/K<\sigma$.
This contradicts Lemma~\ref{lem:nonlinear-weighted-face} by
Lemma~\ref{lem:formal-face-obstruction}, proving the bound on
every $\deg p_i$. The principal monomial itself
has decrease $K$ and order $r$, proving the remaining assertions.
\end{proof}

\begin{theorem}
\label{thm:affine-newton-triangle}
Suppose $\op L$ has polynomial coefficients, order $r\ge2$,
eigenvalues $n$ and monic eigenpolynomials with a recurrence of
fixed length in sufficiently large degrees. After a translation,
for a positive divisor $\sigma$ of $r$ and $K=r/\sigma$,
$$
 a_r(x)=cx^{r-K},\quad
 \op L=\op E+\sum_{i=1}^K\op A_i,\quad
 \op A_i=x^{-i}p_i(\op E),\quad
 \deg p_i\le\sigma i,\quad\deg p_K=r.
$$
No commutation hypothesis is imposed on the $\op A_i$.
\end{theorem}

\begin{proof}
Apply Theorem~\ref{thm:general-newton-triangle} with $d=1$.
For $1\le i\le K$, its degree bound is
$\lfloor1+(r-1)i/K\rfloor=\sigma i$.
\end{proof}

There is also a restriction on the largest degree decrease, without
any commutation assumption. Write
$\op L=\op E+\sum_{i=1}^K\op A_i$, where $\op A_i$ lowers degree
by $i$ and $\op A_K\ne0$. In a monic identity $R(\ad_{\op L})(x)=0$
of degree $N$, the component of degree $1-KN$ is
$(\ad_{\op A_K})^N(x)$, so is zero. Conjugation by
$e^{\op A_K/K}$ therefore gives a finite recurrence.
Proposition~\ref{prop:commuting-recurrence}, applied to this single
component, yields
\refstepcounter{equation}\label{eq:last-decrease-root}
$$
 \op A_K=c\op G^K,\quad \op G=x^{-1}g(\op E),\quad
 g(z)\in\mathbb C[z],\quad c\ne0.
 \leqno\textup{(\theequation)}
$$
In particular, $K$ divides the differential order of $\op A_K$.
The operator $\op G$ itself need not have polynomial coefficients.

\begin{lemma}
\label{lem:poisson-perturbation}
Let $K\ge2$, $\sigma\ge1$, $r=K\sigma$, and
$$
 F_r(z)=\prod_{\ell=1-r}^{1}(z-\ell),\quad
 u=x\xi,\quad v=x^{-1}u^\sigma.
$$
For $c\ne0$, $1\le j<K$ and an integer $a>0$, put
$$
 h_\varepsilon=u+cv^K+\varepsilon v^ju^{-a}.
$$
The coefficient of $\varepsilon$ in
$F_r(\op H_{h_\varepsilon})x$ is nonzero.
\end{lemma}

\begin{proof}
Since $\{u,v\}=-v$, the equations for $\varepsilon=0$ give
$$
 v(t)=v_0z,\quad u(t)=H-cv_0^Kz^K,\quad
 x(t)=\frac{u(t)^\sigma}{v_0z},\quad z=e^{-t},
$$
where $u_0=u(0)$, $v_0=v(0)$ and $H=u_0+cv_0^K$. Choose $u_0,v_0,H\ne0$.
The unperturbed $x(t)$ is a linear combination of the powers
$z^{-1},z^{K-1},\ldots,z^{r-1}$, so
$F_r(\partial_t)x(t)=0$.

Keep the initial point fixed and differentiate the solution
with respect to $\varepsilon$ at zero. Write $u_1,v_1,w$ for
these derivatives of $u,v,x$, respectively.
The equation for $v$ gives
$$
 v_1(t)=a v_0^{j+1}z I(z),\quad
 I'(z)=-\frac{z^{j-1}}{(H-cv_0^Kz^K)^{a+1}},\quad I(1)=0.
$$
Conservation of $h_\varepsilon$ gives
$$
 u_1=u_0^{-a}v_0^j-u^{-a}v^j-cKv^{K-1}v_1.
$$
Consequently
$$
 w=\frac{\sigma u^{\sigma-1}}{v}
       (u_0^{-a}v_0^j-u^{-a}v^j)
 -\frac{u^{\sigma-1}}{v^2}
       \bigl(H+c(r-1)v^K\bigr)v_1.
$$
Thus $w=R(z)+P(z)I(z)$, with $R,P$ rational and $P\ne0$.

If the asserted coefficient were zero, differentiating
$F_r(\partial_t)x(t)=F_r(\op H_{h_\varepsilon})x$ along the solution
at $\varepsilon=0$ would give $F_r(\partial_t)w=0$: the unperturbed
expression is identically zero. Its characteristic roots are the
distinct integers $1-r,\ldots,1$, so $w$ would be a Laurent polynomial
in $z$, and $I=(w-R)/P$ would be rational.
This is impossible: at any root $z_0$ of
$H-cv_0^Kz^K$, the residue of $I'(z)$ is a nonzero multiple of
$$
 \binom{j/K-1}{a}.
$$
Indeed the local change $y=(z/z_0)^K$ reduces its calculation
to the coefficient of $(y-1)^a$ in $y^{j/K-1}$.
This coefficient is nonzero because $0<j/K<1$.
A derivative of a rational function has zero residue at every pole.
\end{proof}

\begin{theorem}
\label{thm:affine-global-normal-form}
Let $\op L$ have polynomial coefficients, order $r\ge2$, and
monic eigenpolynomials with eigenvalues $n$ in every sufficiently
large degree. They have a recurrence of fixed length exactly when,
after a translation of $x$,
$$
 \op L=\op E+\sum_{i\in S}c_i\op G^i,\quad
 \op G=x^{-1}g(\op E),
$$
where $S$ is a finite nonempty set of positive integers,
$c_i\ne0$, $g(z)$ is a nonconstant polynomial, and every occurring
power $\op G^i$ has polynomial coefficients.
Writing $K=\max S$ and $\sigma=\deg g$, one has $r=K\sigma$.
The recurrence order and minimum differential order are $r$, and
the operator is irreducible.
\end{theorem}

\begin{proof}
After the translation in
Theorem~\ref{thm:affine-newton-triangle}, its bounds and
\eqref{eq:last-decrease-root} give
$$
 \op L=\op E+\sum_{i=1}^K\op A_i,\quad
 \deg p_i\le\sigma i,\quad
 \op A_i=x^{-i}p_i(\op E),\quad
 \op A_K=\op G^K,
$$
after rescaling $\op G$. The case $K=1$
already has the required form. Suppose $K\ge2$.
For each $i<K$, choose $c_i$ matching the leading derivative
coefficient when $\ord\op A_i=\sigma i$, and otherwise put $c_i=0$.
Then
$$
 \op B_i=\op A_i-c_i\op G^i
$$
is zero or has order $d_i<\sigma i$. Set $c_K=1$ and
$\op L_0=\op E+\sum_i c_i\op G^i$.

Theorem~\ref{thm:general-order-bound} and
Lemma~\ref{lem:dual-spectral-support} give
$$
 F_r(-\ad_{\op L'})(x)=0,\quad
 F_r(z)=\prod_{h=1-r}^1(z-h).
$$
Transposition gives $F_r(\ad_{\op L})(x)=0$.
Also $F_r(\ad_{\op L_0})(x)=0$: the finite-difference calculation in
the proof of Proposition~\ref{prop:commuting-recurrence} gives
$F_r(\ad_{\op E})(e^{-\ad_{\op Q}}x)=0$ in the degree completion,
where $\op Q=\sum_i c_i\op G^i/i$ and
$\op L_0=e^{\ad_{\op Q}}\op E$.

Suppose some $\op B_i\ne0$. For each nonzero difference put
$\delta_i=\sigma i-d_i$. Choose the smallest $\delta_i$, call it
$\delta$, and let $j$ be the largest index attaining it. Thus $j<K$. Put $d=\sigma j-\delta$ and write the
principal symbols of $\op A_K,\op B_j$ as
$$
 c x^{K(\sigma-1)}\xi^{K\sigma},\quad
 b x^{d-j}\xi^d,\quad cb\ne0.
$$
For a polynomial homogeneous operator $\op C_t$ of degree
change $t$, one has $[\op E',\op C_t']=-t\op C_t'$.
Thus distinct degree changes in the dual identity vanish separately.
The following comparisons take place in $\mathbb C(x)[\partial]$,
without transposing rational components on $\Pp'$.

To isolate the chosen difference $\op B_j$, fix $0\le a\le r$
and consider the terms that change degree by $1-(Ka+j)$ in
$F_r(\ad_{\op L})(x)$. Expand the iterated commutators.
A term of this degree has
at least $a+1$ factors other than $\op E$. If their number is $t$,
its differential order is at most
$$
 \sigma(Ka+j)-t-\sum\delta_i,
$$
where the sum runs over the factors $\op B_i$. Terms containing
none cancel by the identity for $\op L_0$.
Among the others, the largest possible order is
$M_a=\sigma(Ka+j)-(a+1)-\delta$.
Two factors $\op B_i$ give smaller order, as does a larger
$\delta_i$. For one factor with $\delta_i=\delta$, an index below $j$ requires
at least $a+2$ factors to reach total decrease $Ka+j$; an index
above $j$ is excluded by its choice. Thus order $M_a$ receives
only one $\op B_j$, $a$ copies of $\op A_K$, and copies of $\op E$.
At this largest differential order, commutators are calculated by
the Poisson brackets of their leading expressions. The resulting
coefficient is the corresponding term linear in $\varepsilon$ in
$$
 F_r\bigl(\op H_{u+c v^K+\varepsilon b v^j u^{-\delta}}\bigr)
 \left(\frac{u^\sigma}{v}\right),\quad
 u=x\xi,\quad v=x^{-1}u^\sigma.
$$
Lower terms have smaller differential order.
Every linear term uses one perturbation and at most $r$ copies
of $\op A_K$, so the range $0\le a\le r$ includes them all.
Since $d\ge0$, these derivations preserve $\mathbb C(x)[\xi]$;
coefficients with $M_a<0$ are zero.
Lemma~\ref{lem:poisson-perturbation} makes one of the remaining
coefficients nonzero, contradicting $F_r(\ad_{\op L})(x)=0$.
Therefore every $\op B_i$ vanishes. Each occurring $\op G^i$
is polynomial because $c_i\op G^i=\op A_i$.

Conversely, Proposition~\ref{prop:commuting-recurrence} gives the
monic eigenpolynomials, extreme shifts $1,1-r$, and minimum order
$r$ on every high-degree range.
Corollary~\ref{cor:width-irreducibility} gives irreducibility.
\end{proof}

\begin{theorem}
\label{thm:affine-order-equality}
Let $\op L$ have polynomial coefficients, order $r\ge2$ and
nonconstant affine eigenvalues. If its monic eigenpolynomials have a
recurrence of fixed length in sufficiently large degrees, its extreme
shifts are $1,1-r$. Their difference and the minimum nonscalar
differential order are both $r$, even after deleting any finite
initial block. The family is Darboux irreducible.
\end{theorem}

\begin{proof}
Normalize the eigenvalues to $n$.
Theorem~\ref{thm:affine-global-normal-form} and
Proposition~\ref{prop:commuting-recurrence} give the extreme shifts
and minimum order on every high-degree range, and
Corollary~\ref{cor:width-irreducibility} gives irreducibility.
\end{proof}

\begin{proposition}
\label{prop:regular-core-normal-form}
Let $\op G=x^{-1}g(\op E)$, with $g(z)$ a nonconstant polynomial.
For a positive integer $k$, the operator $\op G^k$ has polynomial
coefficients exactly when
$$
 z^{\underline k}\quad\text{divides}\quad
 \prod_{j=0}^{k-1}g(z-j).
$$
If $g(0)=0$, this holds for every $k\ge1$. Otherwise it holds
exactly when $g$ has both a negative integer root and a positive
integer root, and $k\ge a_+-a_-$, where $a_-$ is its largest
negative integer root and $a_+$ its smallest positive integer root.

In the normal form of Theorem~\ref{thm:affine-global-normal-form},
put $K=\max S$, $\sigma=\deg g$, and $r=K\sigma\ge2$.
The largest nonnegative integer root $b_*$ of $g$ exists. Set
$$
 \op G_* =x^{-b_*}\op G x^{b_*}=x^{-1}g(\op E+b_*),\quad
 \op L_* =\op E+\sum_{i\in S}c_i\op G_*^i.
$$
Then $\op G_*$ and $\op L_*$ have polynomial coefficients and
$\ord\op L_*=r$. Let $(P_n(x))_{n\ge0}$ and $(Q_n(x))_{n\ge0}$ be the monic
eigenpolynomials of $\op L,\op L_*$, respectively. They satisfy
$$
 \begin{gathered}
 x^{-b_*}\op L x^{b_*}=\op L_*+b_*\op I,\\[7pt]
 P_{n+b_*}(x)=x^{b_*}Q_n(x),\quad n\ge0.
 \end{gathered}
$$
The last recurrence coefficient of $(P_n(x))_{n\ge0}$ is nonzero
for every $n\ge r-1$ exactly when
$g(m)\ne0$ for every positive integer $m$.
The sequence $(Q_n(x))_{n\ge0}$ has this regularity, and
$$
 \gcd\{P_n(x):n\ge N\}=x^{b_*},\quad N\ge b_*.
$$
\end{proposition}

\begin{proof}
The identity
$$
 \op G^k=x^{-k}\prod_{j=0}^{k-1}g(\op E-j)
$$
gives the divisibility condition by expansion in falling factorials.
Equivalently, each integer $n=0,\ldots,k-1$ must have an integer
root of $g$ in $[n-k+1,n]$. A root at zero covers every index.
Otherwise the indices $n=0$ and $n=k-1$ require roots of both signs.
The nearest roots $a_-,a_+$ cover respectively $n\le a_-+k-1$
and $n\ge a_+$; there is no integer root between them. Coverage
is therefore complete exactly when $a_+\le a_-+k$.

Since at least one positive power occurs in $\op L$, this also
proves existence of $b_*$. The polynomial $g(z+b_*)$ vanishes at
zero. Thus $R(z)=g(z+b_*+1)/(z+1)$ is polynomial and
$\op G_*=R(\op E)\partial$. Conjugation gives the displayed
operator identity and preserves order. Monicity and the eigenvalue
$n+b_*$ then give the identity between the polynomials.

The last recurrence coefficient is
$$
 b_{r-1}(n)=(-c_K)^\sigma[z^\sigma]g(z)
 \prod_{j=0}^{r-2}g(n-j).
$$
For $n\ge r-1$, every argument in this product is positive.
Conversely, a positive integer root $m$ makes it zero at
$n=m+r-2$. This proves the criterion. The polynomial $g(z+b_*)$
has no positive integer root, so the sequence $(Q_n(x))_{n\ge0}$ is regular.

A common zero $a$ of a tail of $(Q_n(x))_{n\ge0}$ would make
$(x-a)^{-1}Q_n(x)$ polynomial for all large $n$, contrary to
Lemma~\ref{lem:dual-pole-polynomiality}.
The polynomial identity above gives the asserted gcd.
\end{proof}

\subsection{Singular cases in composite orders}

The example \eqref{eq:singular-fourth-order-reformulation} extends
to every composite order.

\begin{theorem}
\label{thm:composite-singular-counterexamples}
\label{thm:singular-fourth-counterexample}
Let $k,\sigma\ge2$, $r=k\sigma$, and put
$$
 g(z)=(z-1)(z+1)^{\sigma-1},\quad
 \op G=x^{-1}g(\op E),\quad \op L=\op E+\op G^k.
$$
Then $\op L$ has polynomial coefficients, differential order $r$,
and monic eigenpolynomials
$$
 P_n(x)=\sum_{j=0}^{\lfloor n/k\rfloor}
 \frac{\prod_{a=0}^{kj-1}g(n-a)}{k^j j!}x^{n-kj},\quad n\ge0,
$$
with eigenvalues $n$. Their recurrence has extreme shifts $1,1-r$,
with last coefficient
$$
 b_{r-1}(n)=(-1)^\sigma\prod_{a=0}^{r-2}g(n-a),
$$
which is nonzero for $n\ge r$.
The minimum differential order is $r$ after deleting any finite
initial range, and $\op L$ is irreducible. Affine changes of the
variable and affine normalization of the eigenvalues do not put
$\op L$ in a form with a polynomial generator.
\end{theorem}

\begin{proof}
Since $g(-1)=g(1)=0$ and $k\ge2$,
Proposition~\ref{prop:regular-core-normal-form} shows that $\op G^k$
has polynomial coefficients. Its order is $r=k\sigma$ and its principal
coefficient is $x^{r-k}$.
Proposition~\ref{prop:commuting-recurrence}, with $S=\{k\}$ and $c_k=1$,
gives the displayed eigenpolynomials, recurrence and last coefficient,
and minimum order $r$ on every high-degree range.
The last coefficient is nonzero for $n\ge r$, since the only positive
integer root of $g$ is $1$.
Theorem~\ref{thm:affine-order-equality} gives irreducibility.

Suppose that such a polynomial generator has order $\ell$.
A polynomial $q$ in this generator of degree $r/\ell$ has principal
coefficient of degree $r-r/\ell$. Comparing this with $r-k$ forces
$\ell=\sigma$ and $\deg q=k$.
The only zero of the principal coefficient fixes the translation. After scalar
normalization write the polynomial generator as
$\op H=x^{-1}h(\op E)$, with $h$ monic of degree $\sigma$ and
$h(0)=0$. Distinct powers of $\op H$ have distinct degree decreases,
so only its $k$th power can occur. Comparing the resulting operators gives
$$
 \prod_{a=0}^{k-1}h(z-a)=\prod_{a=0}^{k-1}g(z-a).
$$
For $u=g/h$, divide this identity by the one obtained after
replacing $z$ by $z-1$. The intermediate factors cancel and give
$u(z)=u(z-k)$. A rational function with a nonzero period is constant.
Since $g$ and $h$ are monic, this gives $g=h$, contradicting
$g(0)=-1$ and $h(0)=0$.
\end{proof}

The eigenpolynomials of positive degree are obtained from Horozov's
families with polynomial generators \cite{horozov2016,horozov2018vector}
by multiplication by $x$ and a shift of degree. Put
$$
 \op H=x^{-1}g(\op E+1)=\partial(\op E+2)^{\sigma-1},\quad
 Q_n(x)=e^{\op H^k/k}x^n.
$$
Then
$$
 x^{-1}\op Lx=\op E+\op I+\op H^k,\quad
 P_{n+1}(x)=xQ_n(x),\quad n\ge0.
$$
Both differential orders are $r$.

For $k=\sigma=2$, the recurrence of
\eqref{eq:singular-fourth-order-reformulation} is
$$
 xP_n(x)=P_{n+1}(x)-2n(n^2-1)P_{n-1}(x)
 +n(n-1)^2(n-2)(n-3)(n+1)P_{n-3}(x),
$$
with negative-index terms omitted. The family $(Q_n(x))_{n\ge0}$ is the case $\beta=3$
of \cite[Example~6.3, p.~12]{horozov2016}.

\section{Polynomial sources and admissible transformations}
\label{sec:polynomial-sources}

\subsection{Polynomial sources}

The next lemmas transform the rational source of
Proposition~\ref{prop:general-bispectral-reduction} into a polynomial source
with an $x$-recurrence. Here $^*$ denotes the formal adjoint used in
Theorem~\ref{thm:order-product-identity}, not transposition on functionals.

For $\op A=\sum_{j=0}^r a_j(x)\partial^j$, regular singularity at
zero means that $x^{r-j}a_j(x)/a_r(x)$ is holomorphic at zero for
every $j<r$.

\begin{lemma}
\label{lem:spectral-clearing-intertwiner}
Let $\op A,\op B\in\mathcal D$ have positive order and let
$0\ne\op U\in\mathcal Q$ satisfy $\op A\op U=\op U\op B$.
There is a nonzero $f(z)\in\mathbb C[z]$ such that
$\op U f(\op B)\in\mathcal D$.
\end{lemma}

\begin{proof}
The nonzero right ideal
$I=\{\op C\in\mathcal D:\op U\op C\in\mathcal D\}$
has the form $\op T\mathcal D$. Since $\op B I\subset I$,
$\op B\op T=\op T\op H$ for some $\op H\in\mathcal D$.
The identity $\op T^*\op B^*=\op H^*\op T^*$ makes the local kernel
of $\op T^*$ invariant under $\op B^*$. Section~\ref{sec:preliminaries}
gives $f(\op B^*)=\op V\op T^*$ for some $0\ne f\in\mathbb C[z]$. Hence
$f(\op B)=\op T\op V^*\in I$.
\end{proof}

\begin{lemma}
\label{lem:rational-affine-unit-ladder}
Let $m=k\sigma\ge3$, where $k,\sigma$ are positive integers, and let
$\op H\in\mathbb C(x)[\partial]$ have order $m$ and principal
coefficient $c x^{m-k}$, with $c\ne0$. Suppose that:
\begin{enumerate}[label=(\roman*)]
\item $\op H$ has rational eigenfunctions $\phi_n(x)$, for every
sufficiently large integer $n$, with
$\op H\phi_n(x)=n\phi_n(x)$ and $\deg_\infty\phi_n(x)=n+\ell$ for a fixed
integer $\ell$;
\item for weights $\mathrm{wt}(x)=m-1$ and
$\mathrm{wt}(\partial)=k-m+1$, its leading expression is
$$
 F=x\xi+c\bigl(x^{\sigma-1}\xi^\sigma\bigr)^k;
$$
\item it is regular singular at zero, and every local solution of
$(\op H-s)y=0$ is meromorphic at every finite $a\ne0$, for every
$s\in\mathbb C$;
\item $\mathcal Q$ contains a nonzero $\op Q_0$ satisfying
$[\op H,\op Q_0]=\op Q_0$.
\end{enumerate}
Then a finite chain of rational differential Darboux transformations
and conjugations by nonzero rational functions takes $\op H$ to
$$
 \op H_*+\gamma,\quad
 \op H_*=\op E+q(\op G),\quad
 \op G=x^{-1}g(\op E),\quad
 g(0)=0,\quad\deg g=\sigma,\quad\deg q=k,
$$
where $\gamma$ is an integer and $q,g$ are polynomials, with $q(0)=0$. One may choose $g(j)\ne0$ for every positive integer $j$.
The operators $\op G$ and $\op H_*$ have polynomial coefficients.
For every sufficiently large $n$, the transformation carries
$\phi_n(x)$ to a scalar multiple of the polynomial eigenfunction of
$\op H_*$ of degree $n-\gamma$.
\end{lemma}

\begin{proof}
Apply Lemma~\ref{lem:spectral-clearing-intertwiner} to
$\op Q_0\op H=(\op H-1)\op Q_0$. Then
$\op Q=\op Q_0f(\op H)\in\mathcal D$ satisfies $[\op H,\op Q]=\op Q$.
Put $D=\ord\op Q$. At a common ordinary point, let $\chi_0$ be
the characteristic polynomial of $(\op H|_{\ker\op Q})^\vee$.
Cayley--Hamilton and right division give
$\chi_0(\op H)=\op V\op Q$.

By (i), the highest homogeneous component of $\op H$ is
$\op E+c_0$, with $c_0=-\ell$: a component raising degree would do so
for all but finitely many $\phi_n(x)$, while the degree-preserving component
has value $n$ at $n+\ell$. The unique formal solution with a power
as leading term, called slow, is
$$
 f_s(x)=x^{s-c_0}\sum_{a\ge0}u_a(s)x^{-a},\quad u_0=1.
$$
The $a$th recursion step divides only by the nonzero integer $a$,
so $u_a(s)$ is polynomial in $s$. Put $B=x^{\sigma-1}\xi^\sigma$
and $z=x^{-1/(m-1)}$. In expressions $e^Pz^\rho\sum_{a\ge0}c_a z^a$,
with $P\in z^{-1}\mathbb C[z^{-1}]$ and $\rho,c_a\in\mathbb C$,
the factor $e^P$ is the
exponential part. Solutions with the same exponential part form
one formal branch.
The other $m-1$ solutions, called fast, have distinct nonconstant
exponential parts and leading logarithmic derivatives
$$
 \omega x^{-1+k/(m-1)},\quad
 c\omega^{m-1}=-1.
$$
Since $\gcd(k,m-1)=1$, changing the choice of $z=x^{-1/(m-1)}$
permutes the fast branches transitively. Substitute
$y'/y=z^{m-1-k}\sum_{a\ge0}v_a z^a$ in $\op Hy=sy$.
After multiplication by $z^k$, the initial equation is
$v_0+cv_0^m=0$. At a nonzero root its derivative is $1-m\ne0$,
so all subsequent coefficients are determined. The parameter $s$
first enters $v_k$, with coefficient $-1/(m-1)$.
Terms $a<k$ give the $s$-independent exponential part; $a=k$ gives
the logarithmic exponent, of slope $-1/(m-1)$.
Normalize each formal solution by initial coefficient one after removing
its exponential and power factors. The remaining series coefficients are
polynomial in $s$: the recursions and integration divide only by
$1-m$ and nonzero integers, and the power exponents are affine in $s$.
For every $\mu$ and
$a\ge1$, parameter derivatives of orders $0,\ldots,a-1$ of these
$m$ solutions form a formal basis for $(\op H-\mu)^a$ at infinity:
their highest powers of $\log x$ distinguish derivatives within each part.
For each finite set of parameters and derivative orders, work in a
common differential field extension of $\mathbb C((z))$, generated by
these fundamental systems and their derivatives and with constant field
$\mathbb C$ (a Picard--Vessiot extension). Include in particular the
systems of $\chi_0(\op H)$.
If $\op A\op B$ has a full fundamental system in this field, its
image under $\op B$ lies in $\ker\op A$, of dimension at most
$\ord\op A$. Thus $\dim\ker\op B=\ord\op B$.
In particular, $\dim\ker\op Q=D$.
Enlarge the common field to contain fundamental systems of the
factors and reversed equations before each step, keeping constants $\mathbb C$.
The neighborhoods $z^N\mathbb C[[z]]$, $N\in\mathbb Z$, define the
$z$-adic topology on $\mathbb C((z))$. Powers, logarithms and exponentials
here are formal. No analytic convergence is asserted.
For these formal extensions, see \cite[Chapters~1 and~3]{vanderputsinger2003}.

The highest-order coefficient in $[\op H,\op Q]=\op Q$ gives
$$
 \operatorname{lc}(\op Q)=C x^{(\sigma-1)D/\sigma},
$$
so $\sigma\mid D$.
Write $Q$ for its weighted leading expression. The commutator
identity gives $\{F,Q\}=Q$.
Here
$$
 \mathbb C(x,\xi)=\mathbb C(F,B),\quad
 \{F,B\}=-B,\quad
 x=\frac{(F-cB^k)^\sigma}{B}.
$$
Thus $Q=B^{-1}h(F)$ for a rational function $h$.
Weighted homogeneity makes $h(F)$ a monomial $C F^j$, with integer
$j$. Polynomiality in $\xi$ forces $j\ge\sigma$.
Its $\xi$ degree is $\sigma(kj-1)$ and its weight is $kj-1$.
The term of differential order $D$ has weight $D/\sigma$, so
$kj-1\ge D/\sigma$. The $\xi$ degree of the weighted expression
cannot exceed $D$. Hence
$$
 D=\sigma(kj-1).
$$
The component $\op Q_t$ of highest degree change $t$ satisfies
$t\op Q_t=\op Q_t$, by $[\op H,\op Q]=\op Q$ and the
component $\op E+c_0$ of $\op H$. Hence $t=1$.
Its component $xq_0(\op E)$ has degree $j-\sigma$ in $\op E$,
because the term of $B^{-1}F^j$ raising degree by one is
$x(x\xi)^{j-\sigma}$.
Uniqueness of the slow solution gives
$$
 \op Q f_s(x)=q_0(s-c_0)f_{s+1}(x).
$$
A root of multiplicity $h$ of $q_0(s-c_0)$ contributes precisely the
first $h$ parameter derivatives to the slow kernel, of dimension
$j-\sigma$. The remaining dimension is $(m-1)j$; the root permutation
makes it $j$ on each fast branch.

For a nonzero rational differential operator
$\op A=\sum_i a_i(x)\partial^i$, put
$$
 \nu(\op A)=\min_i\bigl(\operatorname{ord}_0 a_i-i\bigr),\quad
 \alpha(\op A)=\operatorname{ord}_0\operatorname{lc}(\op A)
                  -\operatorname{ord}\op A.
$$
These are the lowest degree change at zero and the change from the
principal coefficient. Both add under products, with $\alpha-\nu\ge0$;
equality means regular singularity. For $\op H$, both equal $-k<0$,
so every nonconstant polynomial in $\op H$ is regular singular at zero.
The product $\chi_0(\op H)=\op V\op Q$ is regular
singular at zero; additivity and nonnegativity of $\alpha-\nu$ give
the same for $\op Q$.
Writing $l=kj-1$, their lowest homogeneous components are
$$
 x^{-k}P(\op E),\quad x^{-l}R(\op E),\quad
 \deg P=m,\quad\deg R=D.
$$
Their indicial polynomials are $P(z)$ and $R(z)$; their roots,
counted with multiplicity, are the exponents at zero. The lowest
commutator vanishes, giving
\refstepcounter{equation}\label{eq:ladder-indicial-balance}
$$
 P(z-l)R(z)=R(z-k)P(z).
 \leqno\textup{(\theequation)}
$$
Choose one representative $\alpha$ of each occurring class modulo
$\mathbb Z$. Let $n_\alpha,D_\alpha$ count the roots of $P,R$ in $\alpha+\mathbb Z$, with multiplicity.
Summing roots in this class in
\eqref{eq:ladder-indicial-balance} gives
$l n_\alpha=kD_\alpha$.
Since $\gcd(k,l)=1$,
$$
 n_\alpha=kr_\alpha,\quad
 D_\alpha=lr_\alpha,\quad
 \sum_\alpha r_\alpha=\sigma.
$$
Only classes with $r_\alpha>0$ occur.

Remove every full right factor $\op H-\mu$ from $\op Q$,
replacing $j$ by its current value after each removal.
The reduced operator still right divides a polynomial in $\op H$,
since $\op Qp(\op H)=p(\op H-1)\op Q$ for every $p\in\mathbb C[z]$.
Rational differential operators preserve exponential parts, which cannot
cancel. Their formal kernels and the images of powers of $\op H-\mu$
on generalized kernels therefore split by these parts. A rational
divisor is unchanged by the root permutation, so it contains all fast
parts or none. A nonzero proper common right divisor with $\op H-\mu$
is thus slow of order $1$ or fast of order $m-1$; a divisor containing both types would
be the full factor.

A slow right factor has a solution with rational logarithmic derivative. For a fast right factor
$\op T$ in $\op H-\mu=\op P\op T$, its left quotient $\op P$
has order $1$; applying $\op T$ to the slow solution gives a nonzero
solution of $\op P$ with rational logarithmic derivative.
Its residues at nonzero finite points are integers, by (iii), and
its exponent at zero belongs to one of the classes $\alpha+\mathbb Z$.
Its power at infinity differs from $\mu-c_0$ by an integer.
Therefore $\mu-c_0\in\alpha+\mathbb Z$.
Multiplying the lowest homogeneous components in
$\op H-\mu=\op P\op T$ adds the exponent multiplicities in each
class modulo integers. The exponent multiplicities in class $\beta$
are $\delta_{\alpha\beta}$ for a right factor of order $1$, and
$n_\beta-\delta_{\alpha\beta}$ for a right factor of order $m-1$.

For the reduced $\op Q$, put $D=\ord\op Q$, $V=\ker\op Q$ and
$A=\op H|_V$ at a common ordinary point. For $U\subset V$, write
$U^\perp=\{\ell\in V^\vee:\langle\ell,U\rangle=0\}$.

Local and formal solution fields have constants $\mathbb C$. For
$\mu\in\mathbb C$ and $h\ge1$, put
$$
 \op D_{\mu,h}=\gcd_{\rm right}(\op Q,(\op H-\mu)^h),
$$
the greatest common right divisor in $\mathcal D$, with principal
coefficient one. In either field its kernel is the
common kernel, and has dimension $\ord\op D_{\mu,h}$: formally,
apply the preceding dimension argument to this right factor of
$(\op H-\mu)^h$, whose fundamental system consists of parameter
derivatives; locally, use the initial-value theorem.
Since $\op Q$ right divides a polynomial in $\op H$, these dimensions
determine all Jordan blocks in both spaces.

On each formal branch, normalized parameter derivatives form one
Jordan chain at $\mu$. Its finite invariant kernel is an initial segment.
The root permutation equates the fast lengths, and reduction excludes
simultaneous slow and fast kernels at one eigenvalue. Thus, for each
eigenvalue $\mu$ of $A$, its generalized eigenspace $V_\mu$ satisfies,
for some $h_\mu\ge1$,
$$
 \begin{gathered}
 N_\mu=(A-\mu)|_{V_\mu},\quad
 K_\mu=\ker N_\mu=\ker\gcd_{\rm right}(\op Q,\op H-\mu),\\[7pt]
 d_\mu=\dim K_\mu\in\{1,m-1\},\quad
 \dim\ker N_\mu^a=d_\mu\min(a,h_\mu),\quad a\ge1.
 \end{gathered}
$$
For $1\le a\le h_\mu$, the map
$$
 N_\mu^{a-1}:\frac{\ker N_\mu^a}{\ker N_\mu^{a-1}}
 \longrightarrow K_\mu
$$
is injective by its kernel and surjective by equality of dimensions.
The rational differential operator $(\op H-\mu)^{a-1}$ inducing this map
commutes with counterclockwise monodromy around zero. Regular
singularity identifies the multiplicity of $\beta+\mathbb Z$ with the
dimension of the generalized monodromy eigenspace for $e^{2\pi i\beta}$
\cite[Chapter~5]{vanderputsinger2003}. Thus every layer has the
exponent multiplicities of $K_\mu$, including resonant exponents and
logarithmic solutions. Within $V_\mu^\vee$, finite duality gives
$$
 (V_\mu/K_\mu)^\vee\simeq K_\mu^\perp
 =(A^\vee-\mu)V_\mu^\vee,
$$
since the last image annihilates $K_\mu$ and has codimension $d_\mu$.

Let $s_\alpha,f_\alpha$ count slow and fast quotients at eigenvalues in
$c_0+\alpha+\mathbb Z$, counting each fast quotient once. The total slow
dimension is $j-\sigma$ and each fast branch has dimension $j$, so
$$
 \sum_\alpha s_\alpha=j-\sigma,\quad
 \sum_\alpha f_\alpha=j.
$$
A slow quotient in class $\alpha$ contributes $\delta_{\alpha\beta}$
in class $\beta$, and a fast quotient contributes
$n_\beta-\delta_{\alpha\beta}$. Additivity gives
$$
 D_\alpha=n_\alpha j+s_\alpha-f_\alpha.
$$
Combining this with $n_\alpha=kr_\alpha$,
$D_\alpha=(kj-1)r_\alpha$ yields
$$
 f_\alpha-s_\alpha=r_\alpha>0.
$$

Choose an eigenvalue $\mu$ occurring in $\ker\op Q$ and put
$$
 \begin{gathered}
 \op T=\gcd_{\rm right}(\op Q,\op H-\mu),\quad
 \op Q=\op Q_1\op T,\quad \op H-\mu=\op P\op T,\\[7pt]
 \widetilde{\op H}=\op T\op P+\mu,\quad
 \widetilde{\op Q}=\op T\op Q_1.
 \end{gathered}
$$
Put $K=\ker\op T$, $\widetilde V=\ker\widetilde{\op Q}$ and
$\widetilde A=\widetilde{\op H}|_{\widetilde V}$. The identities
$$
 \op T\op H=\widetilde{\op H}\op T,\quad
 \op Q_1\widetilde{\op H}=(\op H-1)\op Q_1
$$
and kernel dimensions give
$\op T(V)=\ker\op Q_1$ and $\op Q_1(\widetilde V)=K$.
Write $K^\vee[-1]$ for $K^\vee$ with scalar action $\mu-1$.
Transposition gives the exact sequence
$$
 0\longrightarrow K^\vee[-1]
 \xrightarrow{\op Q_1^\vee}\widetilde V^\vee
 \xrightarrow{\op T^\vee}V^\vee
 \xrightarrow{\mathrm{res}_K}K^\vee\longrightarrow0.
$$
Thus $\operatorname{im}\op T^\vee=K^\perp$: a quotient at $\mu$
is replaced by a submodule at $\mu-1$. Let $\chi_V$ and
$\chi_{\widetilde V}$ be the characteristic polynomials of $A$ and
$\widetilde A$, respectively. For $d_\mu=\dim K\in\{1,m-1\}$, exactness gives
$$
 \chi_{\widetilde V}(z)
 =\chi_V(z)\left(\frac{z-\mu+1}{z-\mu}\right)^{d_\mu}.
$$

The products $\op P\op T$ and $\op T\op P$ have the same principal
coefficient. Cancellation of the nonzero leading symbol of $\op T$
in $\op T\op H=\widetilde{\op H}\op T$ also gives the same weighted
expression $F$ for $\widetilde{\op H}$.
At infinity the highest homogeneous component of $\op T$ is
$x^t t_0(\op E)$, with $t\in\mathbb Z$ and $t_0\ne0$.
Except for the finitely many zeros of $t_0(s-c_0)$, applying
$\op T$ to the slow solution therefore
changes its leading power from $x^{s-c_0}$ to $x^{s-c_0+t}$.
Thus $\widetilde c_0=c_0-t$, and the transformed rational eigenfunctions
have degree $n+\widetilde\ell$ for a fixed integer $\widetilde\ell$
and all large $n$.

Let $y_i(s),\widetilde y_i(s)$ be the normalized formal solutions,
indexed by their exponential parts, $i=1,\ldots,m$. The preceding preservation at
infinity gives polynomial series coefficients in $s$ for both operators.
For $s\ne\mu$, $\op T$ has inverse $(s-\mu)^{-1}\op P$ on solution
spaces and preserves each complete exponential part. Thus the branch
labels agree for all $s$. Write
$$
 \op T y_i(s)=t_i(s)\widetilde y_i(s),\quad
 \op P\widetilde y_i(s)=\tau_i(s)y_i(s),\quad
 t_i(s)\tau_i(s)=s-\mu.
$$
Extraction of the initial normalized coefficient makes $t_i,\tau_i$
holomorphic in $s$. Let $\delta_{ab}$ be the Kronecker delta. Their product gives
$\operatorname{ord}_{s=\lambda}t_i(s)=\varepsilon_i\delta_{\lambda,\mu}$,
where $\varepsilon_i$ is one on the branches of $K$ and zero elsewhere.
Write $\op Qy_i(s)=\kappa_i(s)y_i(s+1)$ and
$\widetilde{\op Q}\widetilde y_i(s)=\widetilde\kappa_i(s)\widetilde y_i(s+1)$.
These scalars are holomorphic by the same coefficient comparison and not
identically zero, since both differential kernels are finite-dimensional.
Differentiating these identities shows that the zero orders $h_i(\lambda)=\operatorname{ord}_{s=\lambda}\kappa_i(s)$ and
$\widetilde h_i(\lambda)=\operatorname{ord}_{s=\lambda}\widetilde\kappa_i(s)$
are the lengths of the Jordan blocks in each branch; zero denotes absence.
The factorizations of $\op Q,\widetilde{\op Q}$ give
$$
 \widetilde\kappa_i(s)
 =\kappa_i(s)\frac{t_i(s+1)}{t_i(s)},\quad
 \widetilde h_i(\lambda)=h_i(\lambda)
 -\varepsilon_i\delta_{\lambda,\mu}
 +\varepsilon_i\delta_{\lambda,\mu-1}.
$$
In particular, for the primary spaces of $\widetilde A$,
$\widetilde V_\mu^\vee\simeq(V_\mu/K)^\vee
=(A^\vee-\mu)V_\mu^\vee$: every block loses one layer.

If $j>\sigma$, choose consecutive eigenvalues of opposite types in a
class with $s_\alpha>0$; one exists since $f_\alpha>s_\alpha$.
Move one layer from the higher value down through the unoccupied values.
At the last reversal, the target is $\nu=\mu-1$: the $d_\mu$ moving branches have
length one and the $m-d_\mu$ old branches retain their positive lengths.
Thus
$$
 \dim\ker(\widetilde A-\nu)
 =\dim\frac{\widetilde V_\nu^\vee}
 {(\widetilde A^\vee-\nu)\widetilde V_\nu^\vee}
 =d_\mu+(m-d_\mu)=m.
$$
Consequently $\ker(\widetilde{\op H}-\nu)
\subseteq\ker\widetilde{\op Q}$, also for repeated roots.

Right division now gives
$\widetilde{\op Q}=\widehat{\op Q}
(\widetilde{\op H}-\nu)$. Indeed, a remainder of order less than
$m$ cannot annihilate a full formal fundamental system of an
order-$m$ equation. Removing this factor subtracts $m$ from the order
of the raising operator and replaces $j$ by $j-1$. It removes one
slow quotient and one fast quotient in the same exponent class, so
$f_\alpha-s_\alpha=r_\alpha$ is unchanged. Repeating the construction
requires exactly $j-\sigma$ such reductions.

At zero, the function $\op A\mapsto\alpha(\op A)-\nu(\op A)$ is nonnegative
and additive under multiplication of nonzero rational differential operators. It vanishes for
$\op P\op T=\op H-\mu$, so it vanishes for each factor and for
$\op T\op P$. Hence $\widetilde{\op H}$ is regular singular there.
Write the lowest homogeneous components of $\op P,\op T$ as
$x^{\nu_P}\pi_P(\op E)$ and $x^{\nu_T}\pi_T(\op E)$, respectively.
Their indicial products in the two orders are
$$
 \pi_P(z+\nu_T)\pi_T(z),\quad
 \pi_T(z+\nu_P)\pi_P(z).
$$
Since $\nu_P,\nu_T\in\mathbb Z$ and
$\nu_P+\nu_T=-k<0$, adding $\mu$ leaves the lowest component
unchanged, and these products have the same root multiplicities in
each class modulo $\mathbb Z$. Thus $n_\alpha$ and $r_\alpha$
are preserved. The same product calculation for
$\op Q_1\op T$ and $\op T\op Q_1$ preserves $D_\alpha$ during
an elementary reversal. Together with the integer change in $c_0$,
this preserves the eigenvalue classes used in the counting argument.

Fix a nonzero finite point $a$ and put $t=x-a$.
Eliminating leading Laurent coefficients gives a meromorphic fundamental
system $y_i=t^{\rho_i}(1+O(t))$ of $(\op H-s)y=0$, with distinct
integers $\rho_i$.
Its Wronskian has order $\sum_i\rho_i-m(m-1)/2$, with nonzero
Vandermonde coefficient. Cramer's rule gives
$\operatorname{ord}_a b_j(x)\ge j-m$ for the coefficients in
$\op H-s=cx^{m-k}(\partial^m+\sum_{j=0}^{m-1}b_j(x)\partial^j)$,
so the equation is regular singular.
The term $-s/(cx^{m-k})$ is holomorphic at $a$, so its indicial
polynomial is independent of $s$. The leading order of a nonzero
meromorphic solution is an indicial root. These roots give a pole
bound independent of $s$. The reverse identity
$\op T\op P=\widetilde{\op H}-\mu$ gives a fixed bound for the transformed
solutions away from $\mu$, including the poles of $\op T$.
At $\mu$, take a fundamental system $y_s$ of
$(\widetilde{\op H}-s)y_s=0$ analytic in $s$. The solutions $\op P y_s$ of
$(\op H-s)\op P y_s=0$ have a uniform pole bound.
Cauchy's coefficient formula preserves this bound under differentiation
in $s$. Differentiating
$\op T\op P y_s=(s-\mu)y_s$ at $s=\mu$ therefore proves meromorphy
of $y_\mu$.
All transformed solutions are therefore meromorphic there.
The kernel of $\op T$ is finite-dimensional, so it annihilates only
finitely many prescribed eigenfunctions. After relabeling
the operators as $\op H,\op Q$, the argument repeats with all four
hypotheses preserved.

At $j=\sigma$, $\ker\op Q$ has no slow solution and dimension
$\sigma$ on each fast branch. The root permutation makes their
characteristic polynomial $\chi$ the same. It has degree $\sigma$
and annihilates $\ker\op Q$, so right division makes
$$
 \op R=\chi(\op H)\op Q^{-1}
$$
a differential operator with rational coefficients of order $\sigma$,
with
$[\op H,\op R]=-\op R$.
Normalize its principal coefficient to $x^{\sigma-1}$.

Suppose first $\sigma\ge2$.
For a term of order $a>1$ in an operator satisfying
$[\op A,\op R]=-\op R$, the highest commutator forces its leading
coefficient to be $C x^{(\sigma-1)a/\sigma}$.
Thus $\sigma\mid a$, and a multiple of
$\op R^{a/\sigma}$ removes it.
Successive removal gives
$$
 \op H=q(\op R)+a_1(x)\partial+c(x),\quad
 \deg q=k,\quad q(0)=0.
$$
The remaining highest commutator reads
$\sigma a_1'-(\sigma-1)a_1/x=1$.
Its rational solution is $a_1=x$.
The slow expansion gives uniqueness of rational eigenfunctions up to
scalars. Thus $\op R$ lowers their degrees by one and $c(x)=O(1)$
at infinity.
With $\zeta'/\zeta=c(x)/x$, put $\widetilde{\op R}=\zeta\op R\zeta^{-1}$.
Conjugation gives
$$
 \zeta\op H\zeta^{-1}=\op E+q(\widetilde{\op R}),\quad
 [\op E,\widetilde{\op R}]=-\widetilde{\op R},\quad
 \widetilde{\op R}=x^{-1}b(\op E),\quad\deg b=\sigma.
$$
Comparing coefficients in $[\op E,\widetilde{\op R}]=-\widetilde{\op R}$ gives this form.

This conjugated operator is ordinary at every $a\ne0$.
Choose an analytic eigenfunction there which does not vanish at $a$.
Its product with $\zeta^{-1}$ is a meromorphic eigenfunction of $\op H$.
Hence $\zeta$ is meromorphic at $a$; its logarithmic derivative has only
a simple pole there, with integer residue.
At zero both operators are regular singular. After division by the
principal coefficient, their coefficients of $\partial^{m-1}$ differ
by $m\zeta'/\zeta$, so $\zeta'/\zeta$ has at most
a simple pole at zero. At infinity it is $O(x^{-1})$.
It follows that $\zeta=x^\beta\rho(x)$ with $\rho\in\mathbb C(x)^\times$.
Consequently
$$
 \rho\op H\rho^{-1}
 =\op E+\beta+q\bigl(x^{-1}b(\op E+\beta)\bigr).
$$
The rational eigenfunctions $\rho\phi_n(x)$ have degree $n-\beta$;
therefore $\beta\in\mathbb Z$.
Ordinariness away from zero makes them Laurent polynomials.
In their lowest Laurent term, the coefficient of the degree decrease
$k$ is a product of $k$ consecutive values of $b(z+\beta)$.
One value vanishes. Choose the largest integer root $v$ of this
polynomial and use
$x^{-v}(\rho\op H\rho^{-1})x^v$.
The resulting $g(z)=b(z+\beta+v)$ satisfies $g(0)=0$ and
$g(j)\ne0$ for all positive integers $j$.
The generator $x^{-1}g(\op E)$ has polynomial coefficients.
Its polynomial eigenfunction in degree $n-\beta-v$ equals the
transformed Laurent eigenfunction, up to a scalar, by the unique
slow expansion.

If $\sigma=1$, write $\op R=\partial+b(x)$.
Subtract $q(\op R)$, with $q(0)=0$, so that the remainder is
$a_1(x)\partial+c(x)$.
The remaining commutator gives
$a_1=x+C$ and $c=(x+C)b+\beta$; hence
$$
 \op H=q(\op R)+(x+C)\op R+\beta.
$$
The degree decrease of $\op R$ gives $b(x)=O(x^{-1})$.
Multiplication by a solution of $\zeta'/\zeta=b$ conjugates this
operator to
$q(\partial)+(x+C)\partial+\beta$.
At each nonzero finite point, the same ratio of local solutions
shows that $\zeta$ is meromorphic.
At zero, $\zeta\phi_n(x)$ solves the conjugated equation, which is ordinary
there. Thus $\zeta\phi_n(x)$ is holomorphic and
$\zeta=(\zeta\phi_n(x))/\phi_n(x)$ is meromorphic.
Hence $\zeta$ is rational. Its transformed eigenfunctions are polynomials
of degree $n-\beta$, so $\beta\in\mathbb Z$. Replacing $q(z)$ by
$q(z)+Cz$ gives the conclusion with $\op G=\partial$.
\end{proof}

The factor construction and recovery of the original degree labels in
\cite[Lemma~3.8, pp.~23--26]{castilloshapiro2026} are adapted below
to Horozov's polynomial sources.

\begin{lemma}
\label{lem:listed-source-exact-bhy}
Let
$$
 \op K=\op E+q(\op G),\quad
 \op G=x^{-1}p(\op E)=\partial h(\op E),\quad
 p(z)=zh(z),
$$
where $q,p\in\mathbb C[z]$ are nonconstant, $q(0)=p(0)=0$,
$(\deg q)(\deg p)\ge2$, and $p$ has no positive integer root.
Suppose that
$\op L$ is a polynomial differential operator with the prescribed
monic eigenpolynomials $P_n(x)$ in every sufficiently large degree,
and
$$
 \op L=\op U\lambda(\op K+\beta)\op U^{-1}
$$
for $0\ne\op U\in\mathcal Q$, $\beta\in\mathbb C$ and
nonconstant $\lambda\in\mathbb C[z]$, where
$\op L P_n(x)=\lambda(n)P_n(x)$.
Then $\beta$ is an integer. A parameter shift of $\op K$ of the
same order gives a polynomial source with a fixed $x$-recurrence
and a finite bispectral Darboux transformation to the
prescribed family, up to nonzero scalars depending on $n$.
\end{lemma}

\begin{proof}
Use Horozov's source polynomials and raising and lowering operators
\cite{horozov2018vector}. Write $q(z)=\sum_{j\ge1}q_jz^j$ and put
$Q(z)=\sum_{j\ge1}q_jz^j/j$ and $\psi_n(x)=e^{Q(\op G)}x^n$, $n\ge0$.
The raising operator is $\op M=e^{\ad_{Q(\op G)}}x$.
The operator $\op G$ lowers the index with coefficient $p(n)$, and
$$
 \op M\op G=p(\op K),\quad
 \op G\op M=p(\op K+1).
$$
Lemma~\ref{lem:spectral-clearing-intertwiner}, applied to
$\op U$ and $\op U^{-1}$, gives nonzero
$\op B,\op C\in\mathcal D$ with
$\op C\op B=c(\op K)$ for a nonzero polynomial $c$.
Put $\Lambda=\lambda(\op K+\beta)$ and
$W=(\op B\Pp+\Pp)/\Pp\subset\mathbb C(x)/\Pp$.
The fixed coefficient poles of $\op B$ make $W$ finite-dimensional;
give it its usual topology. Since $\op L\Pp\subset\Pp$ and
$\op L\op B=\op B\Lambda$, $\op L$ induces an endomorphism
$\overline L$ of $W$. The residue functionals
$[v]\mapsto\operatorname{Res}_{x=a}(x-a)^{j-1}v(x) dx$,
with $a\in\mathbb C$ and $j\ge1$, separate $W$.
Let $\chi$ be the characteristic polynomial of $\overline L^\vee$.
Cayley--Hamilton gives, for $\ell\in W^\vee$,
$$
 \chi\bigl(\lambda(n+\beta)\bigr)\duality{\ell}{[\op B\psi_n(x)]}
 =\duality{\chi(\overline L^\vee)\ell}{[\op B\psi_n(x)]}=0.
$$
The scalar is nonzero for all large $n$, so $\op B\psi_n(x)$ is polynomial.
Its leading term at infinity gives degree $n+b$ for a fixed integer $b$.
Comparison at the highest degree gives
$\lambda(n+b)=\lambda(n+\beta)$, hence $b=\beta$.
Uniqueness now gives a nonzero scalar multiple of $P_{n+\beta}(x)$.
If $\beta=t>0$, replace $(\op B,\op C)$ by
$(\op B\op G^t,\op M^t\op C)$; if $\beta=-t<0$, use
$(\op B\op M^t,\op G^t\op C)$.
The new product $\op C\op B$ belongs to $\mathbb C[\op K]$ by the
raising and lowering identities.
These replacements give maps at the same index:
$$
 \op B\psi_n(x)=a(n)P_n(x),\quad
 \op C P_n(x)=b(n)\psi_n(x).
$$
Their scalar factors are nonzero in all sufficiently large degrees.

The source has a full recurrence with $w=\ord\op K-1\ge1$,
and its last coefficient is nonzero for every $n\ge w$, by
Proposition~\ref{prop:commuting-recurrence} and the absence of positive
integer roots of $p$. Since $\op B\psi_n(x)$ is polynomial in all
sufficiently large degrees, Lemma~\ref{lem:dual-pole-polynomiality}
shows that $\op B$ has polynomial coefficients. Degree comparison in
$\op B\psi_n(x)=a(n)P_n(x)$ shows that $\op B\Pp_j\subseteq\Pp_j$
for every $j$. Let $(\mathbf a_n)_{n\ge0}$ be the source dual basis and
$(\widehat{\mathbf a}_n)_{n\ge0}$ the dual of a monic completion of the target family.

Choose an integer $s\ge\max(1,\ord\op B)$ and put
$$
 \op F_s=\prod_{i=1}^s h(\op E+i),\quad
 \op G_s=\partial h(\op E+s),\quad
 \op K_s=\op E+q(\op G_s).
$$
Direct multiplication gives
$$
 \begin{gathered}
 \op K\op F_s=\op F_s\op K_s,\quad
 \partial^j\op F_s
 =\op G_s^j\prod_{i=1}^{s-j}h(\op E+i),\quad 0\le j\le s.
 \end{gathered}
$$
Thus $\op T=\op B\op F_s$ belongs to
$\mathbb C\langle x,\op E,\op G_s\rangle$.
Let $\psi^{(s)}(x)=(\psi_n^{(s)}(x))_{n\ge0}$ be the monic eigenfamily of $\op K_s$.
Proposition~\ref{prop:commuting-recurrence} gives the partner for $x$;
$\op G_s$ lowers the index and $\op E=\op K_s-q(\op G_s)$.
Hence $x,\op E,\op G_s,\op T$ belong to the source Fourier algebra.
In the weak dual basis $(\mathbf a_j^{(s)})_{j\ge0}$ of $(\psi_j^{(s)}(x))_{j\ge0}$,
$\op K_s'\op F_s'\mathbf a_n=n\op F_s'\mathbf a_n$ gives
$$
 \op F_s'\mathbf a_n=d_s(n)\mathbf a_n^{(s)},\quad
 d_s(n)=\prod_{i=1}^s h(n+i)\ne0,\quad n\ge0.
$$
Pairing with $x^n$ gives the scalar. Dual separation and
$\op T=\op B\op F_s$ give $\op T\psi_n^{(s)}(x)=a_s(n)P_n(x)$,
$a_s(n)=a(n)d_s(n)$, in all large degrees. Since $a(n)=[x^n]\op Bx^n$, $a_s$ is a nonzero
polynomial. Let $\op J_s$ be the source recurrence and $\op b_s$ its Fourier map.
Let $\op M_s$ be their finite raising operator. The identities
$$
 \op R_s=\op M_s^s\partial^s,\quad
 \op R_s\op F_s
 =\prod_{j=0}^{s-1}p_s(\op K_s-j),\quad
 p_s(z)=zh(z+s),
$$
give a reverse differential operator $\op R=\op R_s\op C$ with
rational coefficients and $\op R\op T=c_s(\op K_s)$ for a nonconstant polynomial $c_s$.
Choose a nonconstant polynomial $\theta(x)$ with sufficiently high zeros at
the coefficient poles of $\op R$ and at zero. Then
$\op V=\op R\theta(x)$ has polynomial coefficients, its
$\partial^j$ coefficient being divisible by $x^j$. Therefore
$\op V\in\mathbb C\langle x,\op E\rangle$, in the same
source Fourier algebra. This gives the factorization
$$
 c_s(\op K_s)=\op V\theta(x)^{-1}\op T.
$$
The factors $\op V,\op T,\theta(x)$ satisfy the bispectral Darboux
definition and recover $P_n(x)$ after scalar normalization. The identity
$\op L\op T=\op T\lambda(\op K_s)$ identifies $\op L$ in the eigenoperator algebra of the transformed family.
If the prescribed family has a finite recurrence $\op J_t$ for $x$, put
$$
 \op A=a_s(n)^{-1}\op b_s^{-1}(\op T),\quad
 \op B_\Delta=\op b_s^{-1}(\op V)c_s(n)^{-1}a_s(n).
$$
Reordering as $\theta(x)=\op T c_s(\op K_s)^{-1}\op V$ and using the
diagonal inverse modulo finite rank gives
$\theta(\op J_s)=\op B_\Delta\op A$, as in
Proposition~\ref{prop:general-bispectral-reduction}. Thus
$\op A\psi^{(s)}(x)=P(x)$ and $\op B_\Delta P(x)=\theta(x)\psi^{(s)}(x)$.
The coefficients at shift $i$ in these two factors are, respectively,
$$
 \duality{\mathbf a_{n+i}^{(s)}}{P_n(x)},\quad
 \duality{\widehat{\mathbf a}_{n+i}}{\theta(x)\psi_n^{(s)}(x)}.
$$
Multiplication by $x$ and $\theta(x)$, followed by these coordinate
functionals, gives
$$
 \op J_t\op A=\op A\op J_s,\quad
 \theta(\op J_t)=\op A\op B_\Delta.
$$
Both recurrences represent multiplication by $x$; hence
$\ker\op A$ is $\op J_s$-invariant. All factors are fixed and finite;
the nonzero rational normalizers exclude only finitely many indices.
\end{proof}

\begin{proof}[Proof of Theorem~\ref{thm:all-polynomial-sources}]
Put $d=\deg\lambda$, $r=\ord\op L$ and let $\lambda_d=[z^d]\lambda(z)\ne0$
be the leading coefficient of $\lambda(z)$.
Theorem~\ref{thm:order-product-identity} gives $r=dm$ and a fraction $\op N$ whose class is diagonal with entries $n$, with
$\lambda(\op N)=\op L$.
Proposition~\ref{prop:general-bispectral-reduction} gives a differential operator
$\op H$ with rational coefficients of order $m$, rational eigenfunctions
$\phi_n(x)$ with eigenvalue $n$ and degree at infinity $n+\ell$, and
differential operators $\op B,\op C$ intertwining
$\lambda(\op H)$ and $\op L$ in both directions, with $\op C\op B=f(\op H)$
for a nonzero polynomial $f$.

Theorem~\ref{thm:general-newton-triangle} gives
$a_r=c_0x^{r-K}$ and $K=r/\sigma$.
Since the principal coefficient of $\op H$ is rational and its
$d$th power is $a_r/\lambda_d$, one has $K=dk$,
$m=k\sigma$, and $\operatorname{lc}(\op H)=cx^{m-k}$.
For weights $m-1,k-m+1$, the Newton bound gives
$\mathrm{wt}(\op L)=dk$, with leading expression $F(x,\xi)$.
Canceling the intertwiner's nonzero symbol gives $\mathrm{wt}(\op H)=k$
and leading expression $h$ satisfying
$$
 h^d=\frac{F(x,\xi)}{\lambda_d}.
$$
A rational function whose positive integer power is polynomial is
polynomial. Hence $h\in\mathbb C[x,\xi]$ and $\deg_\xi h\le m$.
A monomial $x^{j-i}\xi^j$ in it has $i\ge0$ and $0\le j\le m$,
since $\op H$ does not raise degree at infinity. It satisfies
$kj-(m-1)i=k$. Since $\gcd(k,m-1)=1$, the only possibilities are
$(i,j)=(0,1),(k,m)$. The degree-zero part of $\op H$ is
$\op E-\ell$, so
$$
 h=x\xi+c x^{m-k}\xi^m.
$$

At zero, the target Newton triangle gives
$\alpha(\op L)=\nu(\op L)=-dk$ in the notation of
Lemma~\ref{lem:rational-affine-unit-ladder}.
Additivity in the intertwining identity gives
$\alpha(\lambda(\op H))=\nu(\lambda(\op H))=-dk$.
Since $\nu(\op H)\le\alpha(\op H)=-k<0$, the highest power of
$\op H$ in $\lambda(\op H)$ determines both valuations. Thus
$\alpha(\op H)=\nu(\op H)=-k$, so $\op H$ is regular singular at zero.
Away from zero, $\op L$ is ordinary. For $s$ outside a finite exceptional set, the
two intertwining operators identify the solutions of $\op H-s$
with a subspace of those of $\op L-\lambda(s)$.
Their poles have a fixed bound, given by the coefficients of $\op C$.
For a fundamental solution $y_s$ analytic in $s$,
$\op C\op B y_s=f(s)y_s$.
If $f$ has a zero of order $h$ at $s_0$, differentiating $h$ times
at $s_0$ gives $f^{(h)}(s_0)y_{s_0}$ on the right. The left has
the same fixed pole bound, so every exceptional solution is also
meromorphic.

Let $S$ be the finite recurrence support, including $1$. In
$\mathcal Q$ put
$$
 \op X=\prod_{a\in S\setminus\{1\}}
 \frac{\ad_{\op N}-a}{1-a}(x).
$$
Modulo finite rank, this is the index-raising term of multiplication
by $x$, so $[\op N,\op X]-\op X$ has zero class. The embedding of
$\mathcal Q$ in this quotient is injective, hence $[\op N,\op X]=\op X$.
Conjugation by the normalization intertwiner gives a nonzero $\op Q_0$
with $[\op H,\op Q_0]=\op Q_0$. All hypotheses of
Lemma~\ref{lem:rational-affine-unit-ladder}
hold; composing its conclusion with the normalization gives
$$
 \op L=\op U\lambda(\op K+\gamma)\op U^{-1},\quad
 0\ne\op U\in\mathcal Q,\quad\gamma\in\mathbb Z,
$$
with $\op K$ of the displayed form.
Lemma~\ref{lem:listed-source-exact-bhy} gives the required finite
factorization and recovers the prescribed $P_n(x)$ at the same indices.
The recurrence order and its nonvanishing coefficient follow from
Theorem~\ref{thm:affine-order-equality} and
Proposition~\ref{prop:regular-core-normal-form}.
\end{proof}

\begin{proof}[Proof of Theorem~\ref{thm:general-irreducible-order-equality}]
Theorem~\ref{thm:general-order-bound} gives $r\ge m$.
For $m=1$, Lemma~\ref{lem:order-one-canonical} gives $r=1$.
For $m=2$, \cite[Theorem~4.1(i), p.~29]{castilloshapiro2026} gives a
polynomial second-order source and a finite chain of recurrence
factor reversals with the original degree labels. Every intermediate
family has a fixed nonscalar polynomial differential eigenoperator,
also in the singular restoration steps
\cite[Lemma~3.9 and proof of Theorem~4.1, pp.~26--28, 33--35]{castilloshapiro2026}.
Thus the chain belongs to
the recurrence-side Darboux class used here.
For $m\ge3$, Theorem~\ref{thm:all-polynomial-sources} gives a
polynomial source of order $m$ with an $x$-recurrence of order $m$.
Both source and target solve the fixed-length problem, and the
transformations preserve the prescribed degree labels. Irreducibility
excludes $r>m$, so $r=m$.
\end{proof}

\begin{proof}[Proof of Theorem~\ref{thm:general-irreducible-classification}]
If the family is Darboux irreducible,
Theorem~\ref{thm:general-irreducible-order-equality} gives recurrence
order $r\ge3$. Theorem~\ref{thm:order-product-identity} then forces
an affine eigenvalue polynomial. Normalize it to $n$ and apply
Theorem~\ref{thm:affine-global-normal-form}.
Conversely, that theorem and
Theorem~\ref{thm:affine-order-equality} give the minimum order,
recurrence and irreducibility for every displayed operator.
\end{proof}

Order one gives translated canonical powers. In order two the Bochner families include Jacobi and Bessel
with quadratic eigenvalues, as well as the affine-spectrum Hermite
and Laguerre cases \cite{bochner1929}.

\subsection{Admissible transformations}

The arguments for removing poles extend
\cite[Lemma~3.6, p.~21, and Lemma~4.3, pp.~30--31]{castilloshapiro2026}.
The sufficient condition uses the finite-dimensional quotient construction
of \cite[Lemma~4.4, p.~31]{castilloshapiro2026}.

For a polynomial source $\op K=\op E+q(\op G)$ of order $m\ge2$,
let $\psi(x)=(\psi_n(x))_{n\ge0}$ be its monic eigenfamily. Write $\op J$ for its
recurrence and $\op b$ for its Fourier map. Put
$$
 \op B=g(n)\op T^{-1}.
$$
The identities
$$
 \op b^{-1}(\op G)=\op B,\quad \op b^{-1}(\op E)=n-q(\op B),\quad \op b^{-1}(x)=\op J
$$
compute the finite difference partner of every element of
$\mathbb C\langle x,\op E,\op G\rangle$.
The inverse Fourier map reverses the order of factors.

\begin{lemma}
\label{lem:admissible-recurrence-test}
Let $s\in\mathbb Z_{\ge0}$ and let $\op F\in\mathbb C\langle\op E,\op G\rangle$ have the form
$$
 \op F=\sum_{j=0}^{s}\op G^j t_j(\op E),\quad
 t_j(z)\in\mathbb C[z],\quad t_0\ne0.
$$
Put $\alpha(n)=t_0(n)$,
$P_n(x)=\alpha(n)^{-1}\op F\psi_n(x)$ and
$\op A=\alpha(n)^{-1}\op b^{-1}(\op F)$.
For sufficiently large $N$, the column $P(x)=(P_n(x))_{n\ge N}$ is monic,
with $\deg P_n(x)=n$. It has a fixed recurrence for multiplication by $x$ precisely when
there is a finite rational difference operator $\widehat{\op J}$ with
\refstepcounter{equation}\label{eq:admissible-recurrence-quotient}
$$
 \widehat{\op J}\op A=\op A\op J.
 \leqno\textup{(\theequation)}
$$
The quotient is unique, has coefficient one at its largest shift,
and has the same extreme shifts as $\op J$. Comparing coefficients reduces its existence to finitely
many rational identities.
\end{lemma}

\begin{proof}
Only the diagonal term $\alpha(n)$ preserves degree, so $\op A$ has
largest shift zero with coefficient one and $P(x)=\op A\psi(x)$.
A recurrence gives $(\widehat{\op J}\op A-\op A\op J)\psi(x)=0$;
independence yields $\widehat{\op J}\op A=\op A\op J$.
Conversely, that identity gives
$xP(x)=\widehat{\op J}P(x)$.
Extreme shifts add under multiplication, so the quotient has largest
shift one and smallest shift $1-m$, with leading coefficient one.
The coefficients at shifts $0,-1,\ldots,1-m$ determine the quotient
successively. The remaining coefficients must vanish; these are the
finite rational identities to be checked.
\end{proof}

\begin{theorem}
\label{thm:polynomial-differential-admissibility}
Let
$$
 \op K=\op E+q(\op G),\quad \op G=x^{-1}g(\op E),\quad
 q(0)=g(0)=0,
$$
where $q,g$ are nonconstant polynomials and $(\deg q)(\deg g)\ge2$.
Let $(\psi_n(x))_{n\ge0}$ be its monic eigenfamily with
eigenvalues $n$. Suppose $\op F$ has polynomial coefficients and does
not increase degree. Define its diagonal polynomial by
$\alpha(n)=[x^n]\op F x^n$ for $n\ge0$, and assume $\alpha\ne0$. Suppose also that
$$
 \op R\op F=c(\op K),\quad
 \op R\in\mathbb C(x)[\partial],\quad
 c\in\mathbb C[z]\setminus\mathbb C.
$$
Put $P_n(x)=\alpha(n)^{-1}\op F\psi_n(x)$ in sufficiently large
degrees. Put
$\beta=\max\{j\in\mathbb Z_{\ge0}:g(-j)=0\}$.
Then the prescribed sequence has a nonscalar polynomial differential
eigenoperator precisely when every coefficient of $\op R$ belongs
to $x^{-\beta}\Pp$.

If this condition holds, put
$$
 \chi_\beta(z)=\prod_{j=1}^{\beta}(z+j),\quad
 \chi_0(z)=1,\quad p(z)=c(z)\chi_\beta(z).
$$
The polynomial differential operator
$\op F\chi_\beta(\op K)\op R$ has eigenvalues $p(n)$ on $P_n(x)$.
\end{theorem}

\begin{proof}
\emph{Necessity.}
Suppose $\op L P_n(x)=\lambda(n)P_n(x)$ for a nonscalar polynomial
differential $\op L$. Its eigenvalue $\lambda(z)$ is a nonconstant polynomial.
The identities on the prescribed polynomials give
$$
 \op L\op F=\op F\lambda(\op K),\quad
 \op R\op L=\lambda(\op K)\op R.
$$
The second follows by right cancellation of $\op F$.
For $v\in\mathbb C(x)$, define its residue functional by
$$
 \duality{\boldsymbol\eta_v}{p(x)}
 =\sum_a\operatorname{Res}_{x=a}v(x)p(x) dx,
 \quad p\in\Pp,
$$
where the sum runs over the finite poles of $v$.
A principal part $c(x-a)^{-j}$ contributes
$c p^{(j-1)}(a)/(j-1)!$. Hermite interpolation gives
$\boldsymbol\eta_v=0$ precisely when $v\in\Pp$.
Put $\op A=\lambda(\op K)$ and
$\Xi=\{\boldsymbol\eta_{\op Rp}:p\in\Pp\}$.
The coefficient poles of $\op R$ bound the supports and orders of
these point functionals, so $\Xi$ is finite-dimensional, with its usual topology.
Formal integration by parts gives
$$
 \boldsymbol\eta_{\op Av}=(\op A^*)'\boldsymbol\eta_v,\quad
 (\op A^*)'\boldsymbol\eta_{\op Rp}=\boldsymbol\eta_{\op R\op Lp}\in\Xi.
$$
A nonconstant multiple $\chi$ of the characteristic polynomial of
$(\op A^*)'|_\Xi$ therefore gives
$\chi(\op A)\op R\Pp\subset\Pp$.
The principal coefficient of $\chi(\op A)$ does not vanish away from
zero, so a pole there of $\op Rp$ would give a nonzero pole of
$\chi(\op A)\op Rp$.
Write $k=\deg q$ and $D=\deg(\chi\circ\lambda)>0$.
If $\op Rp$ has lowest Laurent exponent $l<-\beta$ at zero, its
image under $\chi(\op A)$ has lowest exponent $l-kD$, with nonzero
coefficient proportional to $\prod_{j=0}^{kD-1}g(l-j)$, since no
integer root of $g$ lies below $-\beta$. This is again a contradiction.
Thus $x^\beta\op R$ maps $\Pp$ to $\Pp$ and has polynomial
coefficients by the coefficient criterion.

\emph{Sufficiency.}
Assume the stated coefficient condition. The space
$x^{-\beta}\Pp$ is $\op G$-invariant because $g(-\beta)=0$.
On $\mathcal M_\beta=x^{-\beta}\Pp/\Pp$,
the operator $\op K$ is triangular with diagonal
entries $-\beta,-\beta+1,\ldots,-1$. Its transpose on the finite dual
has characteristic polynomial $\chi_\beta$.
Cayley--Hamilton and separation by $\mathcal M_\beta^\vee$ give
$\chi_\beta(\op K)\mathcal M_\beta=0$.
The coefficient criterion makes $\chi_\beta(\op K)\op R$ and
$\op F\chi_\beta(\op K)\op R$ polynomial.
Finally,
$$
 \bigl(\op F\chi_\beta(\op K)\op R\bigr)\op F
 =\op F p(\op K),
$$
which proves the eigenvalue assertion. Its order is
$(\ord\op K)\deg p>0$.
\end{proof}

\begin{corollary}
\label{cor:finite-spectral-jet-algebra}
Under the admissibility condition of
Theorem~\ref{thm:polynomial-differential-admissibility}, put
$p=c\chi_\beta$ and $s=\deg p$. Then
$$
 p\mathbb C[z]\subset\mathcal A_F,\quad
 \mathcal A_F=p\mathbb C[z]\mathbin\oplus V,
 \quad V\subset\mathbb C[z]_{<s}.
$$
Here $\mathbb C[z]_{<s}=\{r\in\mathbb C[z]:\deg r<s\}$.
The space $V$ is the kernel of a finite linear system computed from
$\op F,\op K,p$. Modulo $p$, it is a unital subalgebra of
$\mathbb C[z]/(p)$. In particular, the minimum differential order
is computable using only spectral polynomials of degree at most $s$.
The greatest common divisor of all positive differential orders is
$\ord\op K$. The fraction field of $\mathcal A_F$ is $\mathbb C(z)$,
and its integral closure is $\mathbb C[z]$.
\end{corollary}

\begin{proof}
For every $h\in\mathbb C[z]$, the operator
$\op F h(\op K)\chi_\beta(\op K)\op R$ is polynomial and corresponds
to $ph$. Divide $\lambda=ph+r$, with $\deg r<s$.
Then $\lambda\in\mathcal A_F$ precisely when $r\in\mathcal A_F$.
Right division of $\op F r(\op K)$ by $\op F$ over
$\mathbb C(x)[\partial]$ is linear in the coefficients of $r$.
Clearing a fixed denominator, zero remainder and polynomial quotient
give a finite homogeneous linear system defining $V$.
Every polynomial differential eigenoperator has a polynomial eigenvalue
and the same intertwining identity, so no others occur. Composition
gives the algebra assertion.
A residue class in $\mathbb C[z]/(p)$ is determined by the values
and derivatives at the roots of $p$, up to one less than each root's
multiplicity. Thus the finite equations can be written in terms of
these values and derivatives of $\lambda$.
The polynomial $p$ itself is admissible, giving the asserted finite
bound for the search for the smallest positive degree.
Every positive differential order is $(\ord\op K)\deg\lambda$;
the polynomials $p$ and $zp$ give consecutive spectral degrees.
Their orders therefore have greatest common divisor $\ord\op K$.
Since $p,zp\in\mathcal A_F$, its fraction field contains $z$.
The element $z$ is integral over $\mathbb C[p]\subset\mathcal A_F$.
Since $\mathbb C[z]$ is integrally closed, it is the integral closure
of $\mathcal A_F$.
\end{proof}

Rank preservation is a classical Darboux property \cite[Proposition~1]{bhy1996grassmann}.
Iliev determined the full Krall--Laguerre spectral algebra for generic
admissible parameters \cite[Definition~5.1 and Theorem~5.3]{iliev2013}. Here the finite
system determines the full algebra for every factorization satisfying
Theorem~\ref{thm:polynomial-differential-admissibility}, including repeated roots.

\begin{proof}[Proof of Theorem~\ref{thm:complete-admissible-data}]
Assume (i) and (ii), and put $\op R=x^{-\beta}\op S$.
Then $\op R\op F=c(\op K)$, and
Theorem~\ref{thm:polynomial-differential-admissibility} gives the
polynomial differential eigenoperator
$$
 \op L_*=\op F\chi_\beta(\op K)x^{-\beta}\op S,\quad
 \lambda_*(z)=c(z)\chi_\beta(z).
$$
Condition (ii) gives the finite recurrence of order $m$.

Choose $M\ge\max(1,\beta+\ord\op R)$ and put
$\op V=\op R x^M$. Its $\partial^j$ coefficient is polynomial
and divisible by $x^j$: in the product rule, every term has
valuation at least $M-\beta-\ord\op R+j\ge j$ at zero.
Thus $\op V\in\mathbb C\langle x,\op E\rangle$ and
$$
 c(\op K)=\op V x^{-M}\op F
$$
is the required bispectral factorization: $\op V$, $\op F$ and
the multiplier $x^M$ belong to the source Fourier algebra.
With
$$
 \op B_\Delta=\op b^{-1}(\op V)c(n)^{-1}\alpha(n),
$$
its finite difference factors satisfy
$$
 \begin{gathered}
 \op J^M=\op B_\Delta\op A,\quad
 \widehat{\op J}^{M}=\op A\op B_\Delta,\\[7pt]
 \widehat{\op J}\op A=\op A\op J.
 \end{gathered}
$$
The first identity follows from
$x^M=\op F c(\op K)^{-1}\op V$ modulo finite rank, as in
Proposition~\ref{prop:general-bispectral-reduction}; the last is the
recurrence intertwining identity, and right cancellation of $\op A$
gives the second.
After division by $\alpha(n)$, the transformed family is exactly
$P_n(x)$ with its original degree labels.

Suppose first that the displayed family has a nonscalar
polynomial differential eigenoperator $\op L$ and a recurrence of
fixed length. Write $\lambda(z)$ for its eigenvalue polynomial. Then
$\op L\op F=\op F\lambda(\op K)$.
Apply Lemma~\ref{lem:spectral-clearing-intertwiner} to
$\lambda(\op K)\op F^{-1}=\op F^{-1}\op L$.
It gives $\op R=\op F^{-1}f(\op L)\in\mathcal D$, where $f(z)$
may be chosen nonconstant by multiplying it by $z$. Thus
$\op R\op F=c(\op K)$, with $c=f\circ\lambda$ nonconstant.
Theorem~\ref{thm:polynomial-differential-admissibility} gives
$\op R=x^{-\beta}\op S$ with $\op S\in\Pp[\partial]$,
proving (i). Its recurrence gives (ii) by
Lemma~\ref{lem:admissible-recurrence-test}.

Let $(P_n(x))_{n\ge n_0}$ be any solution of recurrence order $m$.
Theorem~\ref{thm:all-polynomial-sources} and
Lemma~\ref{lem:listed-source-exact-bhy} provide a polynomial source
$\op K$ of order $m$, a polynomial forward factor
$\op F\in\mathbb C\langle x,\op E,\op G\rangle$, and a rational
reverse $\op R$ with $\op R\op F=c(\op K)$.
Since $\op F\psi_n(x)$ is a nonzero multiple of $P_n(x)$ for large $n$,
$\op F$ preserves degree. The relations
$$
 x\op G=g(\op E),\quad \op Gx=g(\op E+1),\quad
 [\op E,\op G]=-\op G
$$
reduce every word to terms $x^j a_j(\op E)$, $j\ge0$, or
$\op G^j t_j(\op E)$, $j\ge1$.
The degree changes are distinct. Degree preservation removes the
positive changes and gives the displayed form of $\op F$.
Theorem~\ref{thm:polynomial-differential-admissibility} forces
$\op R=x^{-\beta}\op S$ with polynomial $\op S$, proving (i).
The prescribed recurrence gives (ii) by
Lemma~\ref{lem:admissible-recurrence-test}. The intertwining identity
identifies every original eigenoperator with $\op L_\lambda$ and gives
$\ord\op L_\lambda=m\deg\lambda$ for nonconstant $\lambda$.
Corollary~\ref{cor:finite-spectral-jet-algebra} gives the full algebra
and its minimum order.
\end{proof}

For a fixed source and each fixed $h=\deg c$, the factorization forces
$\ord\op F\le mh$ and $\ord\op S=mh-\ord\op F$.
In the displayed form of $\op F$ one has
$j\deg g+\deg t_j\le mh$.
The $\partial^j$ coefficient of $\op S$ has degree at most
$\beta+j$, because the largest degree change of
$\op S\op F=x^\beta c(\op K)$ is $\beta$.
Thus for fixed source and $h$, (i) is a finite system of polynomial
coefficient equations and (ii) a finite system of rational identities.
The minimum eigenoperator uses spectral degree at most $\deg c+\beta$.

\subsection{The families}

For each divisor $\sigma$ of $m\ge3$, put $k=m/\sigma$.
Rescaling $q$ and the $t_j$ allows $g$ to be monic, so the source
parameters are
$$
 g(z)=z\prod_{\nu=1}^{\sigma-1}(z-\rho_\nu),\quad
 q(z)=\sum_{j=1}^k q_jz^j,\quad q_k\ne0.
$$
The roots $\rho_\nu$ are arbitrary, with repetitions allowed.
For $\sigma=1$ the sources are Appell; for $k=1$ there is one
homogeneous lowering term. Intermediate divisors allow several
powers of the same higher-order $\op G$.

For the admissible data of Theorem~\ref{thm:complete-admissible-data}, set
$$
 \gamma_\ell=[z^\ell]\exp\left(\sum_{j=1}^k\frac{q_j}{j}z^j\right),
 \quad\ell\ge0.
$$
Expanding the exponential and applying the factors in $\op F$ gives
\refstepcounter{equation}\label{eq:classified-family-coefficients}
$$
 \begin{aligned}
 P_n(x)&=\frac1{t_0(n)}\sum_{\ell=0}^n\Gamma_\ell(n)
 \left(\prod_{i=0}^{\ell-1}g(n-i)\right)x^{n-\ell},\\[7pt]
 \Gamma_\ell(n)&=\sum_{j=0}^{\min(s,\ell)}
 \gamma_{\ell-j}t_j(n-\ell+j).
 \end{aligned}
 \leqno\textup{(\theequation)}
$$
Empty products are one and $t_0(n)\ne0$ for large $n$.
For $\op F=\op I$, $\Gamma_\ell=\gamma_\ell$ and $P_n(x)=\psi_n(x)$.

\begin{corollary}
\label{cor:constant-principal-appell-classification}
Let a solution have recurrence order $m\ge3$ and constant principal
differential coefficient. Its polynomials in sufficiently large degrees
are Appell:
$$
 P_n(x)=e^{Q(\partial)}x^n,\quad
 Q(z)=\sum_{j=1}^m\frac{q_j}{j}z^j,\quad q_m\ne0.
$$
Every polynomial differential eigenoperator for these degrees is a
polynomial in $\op E+\partial Q'(\partial)$. The minimum differential
order for these degrees is $m$.
\end{corollary}

\begin{proof}
Translation changes only the linear term of $Q$, so use the
normalization of Theorem~\ref{thm:all-polynomial-sources}.
Choose the source $\op K=\op E+q(\op G)$, with
$\op G=x^{-1}g(\op E)$, and factors $\op F,\op R$ from
Theorems~\ref{thm:all-polynomial-sources} and~\ref{thm:complete-admissible-data}. They satisfy
$\op R\op F=c(\op K)$, with $c\in\mathbb C[z]$ nonconstant,
and $\op F$ has polynomial coefficients and preserves degree.
For the monic source eigenpolynomials $\psi_n(x)$, one has
$$
 P_n(x)=\alpha(n)^{-1}\op F\psi_n(x),\quad
 \alpha(n)=[x^n]\op F x^n.
$$
Let $\lambda(n)$ be the eigenvalue of the given operator $\op L$.
Principal coefficients in
$\op L\op F=\op F\lambda(\op K)$ make the principal coefficient
of $\op K$ constant. Thus $\deg g=1$; rescale to put
$\op G=\partial$. Write
$$
 \op K=\op E+\partial Q'(\partial),\quad
 \psi_n(x)=e^{Q(\partial)}x^n.
$$
Since $g(z)=z$ has no negative integer root,
Theorem~\ref{thm:polynomial-differential-admissibility} makes $\op R$
polynomial.
For an Appell source, every polynomial differential operator has a
finite difference partner: conjugation by $e^{Q(\partial)}$ fixes $\partial$ and
sends $x$ to $x+Q'(\partial)$.
Applying the inverse Fourier map to $\op R\op F=c(\op K)$ gives
$\op b^{-1}(\op F)\op b^{-1}(\op R)=c(n)$.
Extreme shift widths add, so $\op b^{-1}(\op F)$ is a single shift.
Its degree change is zero; hence $\op b^{-1}(\op F)=\alpha(n)$ and
$P_n(x)=\psi_n(x)$.
The spectral algebra follows from
Proposition~\ref{prop:affine-spectral-algebra}.
\end{proof}

For Horozov's sources \cite{horozov2018vector}, the first-order factor
$\op E$ gives the following full eigenoperator algebra.

\begin{proposition}
\label{prop:cubic-admissibility-example}
Let $h(z)$ be a nonzero polynomial and put $m=\deg h+1$. Set
$$
 \op G=x^{-1}\op E h(\op E),\quad
 \op K=\op E+\op G,\quad \psi_n(x)=e^{\op G}x^n,\quad n\ge0,
$$
and define
$$
 P_0(x)=1,\quad P_n(x)=\frac1n\op E\psi_n(x),\quad n\ge1.
$$
This sequence has recurrence order $m$. Its polynomial differential
eigenoperators are exactly $\lambda(\op H)$, where
$$
 \op H=\op E+x^{-1}(\op E-1)h(\op E),\quad
 \lambda\in\mathcal A_{\op E}.
$$
If $h(0)\ne0$ and $h$ has a negative integer root, set
$$
 b=\min\{j\in\mathbb Z_{\ge1}:h(-j)=0\},\quad
 p_b(z)=\prod_{j=0}^{b}(z+j).
$$
Then
$$
 \mathcal A_{\op E}=
 \begin{cases}
 \mathbb C[z],&h(0)=0,\\[7pt]
 \mathbb C+p_b(z)\mathbb C[z],&h(0)\ne0\text{ and }b\text{ exists},\\[7pt]
 \mathbb C,&h(0)\ne0\text{ and }h(-j)\ne0\text{ for all integers }j\ge1.
 \end{cases}
$$
The minimum differential order is $m$ in the first case and
$m(b+1)$ in the second. In the third case there is no nonscalar
polynomial differential eigenoperator. The algebra and these minima
are unchanged if only sufficiently large degrees are prescribed.
\end{proposition}

\begin{proof}
Put $\op G_+=x^{-1}\op E h(\op E+1)$. The identity
$\partial\op G=\op G_+\partial$ gives
$$
 P_n(x)=x e^{\op G_+}x^{n-1},\quad n\ge1.
$$
Proposition~\ref{prop:commuting-recurrence} gives its recurrence
order $m$. The recurrence starts at $n=1$, with no $P_0(x)$ term;
at $n=0$, use $xP_0(x)=P_1(x)$.
The identity $\op H\op E=\op E\op K$ gives $\op H P_n(x)=nP_n(x)$
for $n\ge1$. By Lemma~\ref{prop:unrestricted-eigenpolynomials}, a polynomial
eigenoperator $\op L$ has a polynomial eigenvalue $\lambda$.
The rational differential operator $\op L-\lambda(\op H)$ annihilates a tail
and is zero.

If $h(0)=0$, $\op H$ has polynomial coefficients.
Otherwise, when $b$ exists, both $x\Pp$ and
$$
 W_b=\spanC\{x^{-b},\ldots,x^{-1},1\}
$$
are $\op H$-invariant and $W_b\cap\Pp=\mathbb C1$.
Its bidiagonal action makes $1$ cyclic with annihilator polynomial $p_b$.
Writing $1^\perp=\{\ell\in W_b^\vee:\ell(1)=0\}$, dual separation shows that
$\lambda(\op H)\Pp\subset\Pp$ if and only if
$$
 \duality{\ell}{\lambda(\op H)1}=0,\quad \ell\in1^\perp.
$$
By cyclicity, this is equivalent to $\lambda\in\mathbb C+p_b\mathbb C[z]$.

If $b$ does not exist, a nonconstant $\lambda$ of degree $d$ gives
$$
 [x^{-d}]\lambda(\op H)1
 =[z^d]\lambda(z)\prod_{j=0}^{d-1}(-j-1)h(-j)\ne0.
$$
Thus only constants remain. For polynomial $\lambda(\op H)$, its value
on $P_0(x)=1$ is $\lambda(0)$, so the algebra also applies to the full
sequence. Finally
$\ord\lambda(\op H)=m\deg\lambda$ for nonconstant $\lambda$,
which gives the stated minima.
\end{proof}

For example, take $h(z)=(z+b)^2$, with any integer $b\ge1$.
Then the first-order factor $\op E$ gives
$$
 \mathcal A_{\op E}=\mathbb C+z(z+1)\cdots(z+b)\mathbb C[z].
$$
The minimum differential order is $3(b+1)$, while the recurrence
order is three. Explicitly,
$$
 \begin{aligned}
 \op H_b&=x^2\partial^3+2(b+1)x\partial^2+(x+b^2)\partial-b^2/x,\\[7pt]
 \op L_b&=\prod_{j=0}^b(\op H_b+j\op I),\quad
 \op L_b P_n(x)=n(n+1)\cdots(n+b)P_n(x).
 \end{aligned}
$$
The product $\op L_b$ is polynomial. Arbitrarily large gaps between rank and
minimum order occur in
\cite[Section~5.1, p.~365]{bhy1997}; here they occur with a fixed
four-term recurrence and a first-order factor. \section{Further classifications}
\label{sec:further-classifications}

\subsection{Sheffer families with polynomial eigenvalues}

A monic Sheffer family has a formal generating function
$$
 F(x,z)=\sum_{n\ge0}P_n(x)\frac{z^n}{n!}
 =A(z)e^{x\psi(z)},\quad A(0)=\psi'(0)=1,\quad\psi(0)=0.
$$
Define power series $f(t),g(t)$ by
$$
 f(\psi(z))=z\psi'(z),\quad
 g(\psi(z))=zA'(z)/A(z).
$$
Then $f(0)=g(0)=0$ and $f'(0)=1$. The generator
$\op B=xf(\partial)+g(\partial)$ acts on polynomials by finite sums,
with $\op BP_n(x)=nP_n(x)$ \cite[Theorem~2.3]{sheffer1939}.

\begin{theorem}
\label{thm:sheffer-polynomial-generator}
A monic Sheffer family with a fixed-length recurrence in sufficiently
large degrees has a nonscalar polynomial differential eigenoperator
on those degrees if and only if its generator has one of the forms
$$
 \begin{aligned}
 &\op B=\op E+g(\partial),\quad g(t)\in t\mathbb C[t],\\[7pt]
 &\op B=x(\partial+a\partial^2)+b\partial+c\partial^2,
 \quad a\ne0,\quad b,c\in\mathbb C.
 \end{aligned}
$$
In either case the rational differential eigenoperator algebra on
sufficiently large degrees is $\mathbb C[\op B]$; nonlinear spectra
therefore give no further Sheffer families.
\end{theorem}

\begin{proof}
Put $\pi=1/\psi'$ and $\rho=A'/(A\psi')$.
The identity $xF=\pi\partial_zF-\rho F$ gives the recurrence coefficient
of $P_{n-s}(x)$ as
$$
 n^{\underline s}\bigl((n-s)[z^{s+1}]\pi-[z^s]\rho\bigr).
$$
For every forbidden shift, both coefficients of this affine function
of $n$ vanish. Thus $\pi,\rho$ are polynomials.

Under the moment identification
$\mathbf u\mapsto\sum_{n\ge0}\duality{\mathbf u}{x^n}t^n/n!$,
the weak topology $\sigma(\Pp',\Pp)$ is the product topology of the
coefficients, with their usual complex topologies, and
$x'=\partial_t$, $\partial'=t$.
Let $\op M$ be a polynomial differential eigenoperator with nonconstant
eigenvalue polynomial $\lambda$. Then
$\op A=\op M-\lambda(\op B)=\sum_{j=0}^D x^j h_j(\partial)$,
$h_j\in\mathbb C[[t]]$, annihilates all large $P_n(x)$, so
$\op A'=\sum_{j=0}^D h_j(t)\partial_t^j$ has finite rank.
A nonzero differential operator over $\mathbb C((t))$ has
finite-dimensional kernel, hence infinite rank on $\mathbb C[[t]]$.
Thus $\op A'=0$ and $\op M=\lambda(\op B)$.
Put $d=\deg\lambda\ge1$.

Its coefficient of $x^d$ is $[t^d]\lambda(t)f(\partial)^d$, so
$f(t)^d=Q(t)$ for a nonconstant polynomial $Q$. Since
$$
 Q(\psi(z))=(z\psi'(z))^d=(z/\pi(z))^d,
$$
$\psi$ is algebraic and $\psi'=1/\pi$ is rational, with $\pi(0)=1$.
Integration gives a rational function plus logarithms $c_a\log(z-a)$.
A nonzero $c_a$ would give infinitely many branches, contrary to
algebraicity. Hence $\psi$ is rational.

If $\pi$ is constant, $\psi(z)=z$ and $f(t)=t$.
Otherwise $\psi'\to0$ at infinity, so $\psi=C_0+R$ with $R(\infty)=0$,
$R(z)\sim C_1z^{-m}$. Let its finite pole orders be $q_1,\ldots,q_s$
and put $S=\sum_iq_i$. The derivative has pole orders $q_i+1$,
no finite zeros, and behaves as $-mC_1z^{-m-1}$ at infinity. Hence
$$
 S+s=\deg\pi=m+1\le S+1.
$$
Indeed $R$ has denominator degree $S$ and nonzero numerator, so
$m\le S$. Thus $s=1$, and normalization gives, for $a\ne0$, $q\in\mathbb Z_{\ge1}$,
$$
 \psi(z)=\frac{(1-az)^{-q}-1}{aq}.
$$
For $w=1-az$, $Q(\psi(z))$ is invariant under
$w\mapsto\omega w$, $\omega^q=1$, whereas
$$
 (z\psi'(z))^d=a^{-d}(1-w)^d w^{-d(q+1)}
$$
has its only nonzero finite zero at $w=1$. Hence $q=1$, so
$\psi(z)=z/(1-az)$ and $f(t)=t(1+at)$.

For $f(t)=t$, $\rho$ is polynomial, so $g(t)=t\rho(t)$.
For $f(t)=t(1+at)$ the same condition gives
$$
 g(t)=t(1+at)\rho\left(\frac{t}{1+at}\right),
$$
whose only possible finite pole is $t_0=-1/a$.
Suppose its order is $h\ge1$. On the moment dual,
$\op B'=f(t)\partial_t+g(t)=\op C$.
Since $f$ has a simple zero at $t_0$, $f\partial_t$ does not increase
pole order, whereas multiplication by $g$ increases it by $h$.
Thus $\op C^d1$ has pole order $dh$, with leading term from $g^d$;
smaller powers have smaller pole order. Hence $\lambda(\op C)1$ has a pole.

But for $\op M=\sum_jm_j(x)\partial^j$, the moment series of
$\op M'\boldsymbol\delta_0$ is the polynomial $\sum_jm_j(0)t^j$.
As $\boldsymbol\delta_0$ corresponds to $1$, this contradicts
$\op M'=\lambda(\op C)$. Thus $g$ is polynomial, and the displayed
formula gives $\deg g\le2$.

Conversely, the first form has $\pi=1$, $\rho(z)=g(z)/z$;
the second has $\pi(z)=(1-az)^2$, $\rho(z)=b+(c-ab)z$.
These are polynomials, so the coefficient formula gives a finite
recurrence. Both generators have eigenvalues $n$; apply
Proposition~\ref{prop:affine-spectral-algebra} for the algebra assertion.
\end{proof}

\subsection{One homogeneous lowering term}
\label{sec:single-decrease}

\begin{lemma}
\label{lem:finite-recurrence-poles}
Let $\op M$ have polynomial coefficients and a nonconstant polynomial
eigenvalue $\lambda(n)$, with monic eigenpolynomials in sufficiently
large degrees. Write
$$
 P_n(x)=\sum_{j=0}^n c_j(n)x^{n-j},\quad c_0(z)=1.
$$
Each $c_j$ is rational. A fixed-length recurrence implies that all
their poles lie in $\Omega+\mathbb Z$ for some finite
$\Omega\subset\mathbb C$, and hence have uniformly bounded imaginary parts.
\end{lemma}

\begin{proof}
Lemma~\ref{prop:unrestricted-eigenpolynomials} and coefficient
comparison give rational $c_j,b_s$ and, for $j\ge1$,
$$
 c_j(z)-c_j(z+1)=\sum_{s=0}^{w}b_s(z)c_{j-1-s}(z-s),
$$
where $c_j=0$ for $j<0$. Let $\Omega$ contain the poles of the $b_s$.
Inductively, a pole of $c_j$ outside $\Omega+\mathbb Z$ with maximal
real part in its integer class cannot cancel with $c_j(z+1)$;
the right-hand side is regular there, a contradiction.
\end{proof}

\begin{theorem}
\label{prop:homogeneous-recurrence}
Let $k,\ell$ be integers with $1\le k\le\ell$, and put
$$
 \op{A}=\sum_{s=k}^{\ell}\alpha_sx^{s-k}\partial^s
        =x^{-k}p(\op{E}),\quad
 p(z)=\sum_{s=k}^{\ell}\alpha_sz^{\underline s},\quad\alpha_\ell\ne0.
$$
For nonconstant $\lambda\in\mathbb C[z]$, the monic eigenpolynomials
of $\op M=\lambda(\op E)+\op A$ in sufficiently large degrees have a
fixed-length recurrence exactly in the following cases:
\begin{enumerate}[label=(\roman*)]
\item
\refstepcounter{equation}\label{eq:homogeneous-root-criterion}
$$
 \lambda(z)=az+b,\quad a\ne0,\quad
 p(z)=c\prod_{i=0}^{k-1}g(z-i),
 \leqno\textup{(\theequation)}
$$
where $g(z)$ is a monic nonconstant polynomial and $c\ne0$.
\item $\deg\lambda=2$ and $\ell\le2$.
\end{enumerate}
In case~(i), $\ell=k\sigma$, with $\sigma=\deg g$, and
$P_n(x)=e^{\op{A}/(ak)}x^n$. With
$h(z)=c\prod_{i=1}^{k-1}g(z-i)$ and
$\nabla_kg(z)=g(z)-g(z-k)$, the recurrence is
\refstepcounter{equation}\label{eq:homogeneous-recurrence-coefficients}
$$
 \begin{gathered}
 xP_n(x)=P_{n+1}(x)+\sum_{q=1}^{\sigma}\gamma_q(n)P_{n+1-qk}(x),\\[7pt]
 \gamma_q(n)=\frac{(-1)^qh(n+1)\nabla_k^qg(n+1)}{q!(ak)^q}
 \prod_{i=1}^{q-1}p(n+1-ik),\quad 1\le q\le\sigma.
 \end{gathered}
 \leqno\textup{(\theequation)}
$$
Terms with negative indices are omitted. Whenever $\lambda$ is
affine, $(P_n(x))_{n\ge N}$ has minimum differential order $\ell$
for every sufficiently large $N$. In case~(ii), the recurrence has
three terms.
\end{theorem}

\begin{proof}
Put $d=\deg\lambda$. Coefficient comparison gives
$$
 P_n(x)=\sum_{q\ge0}c_q(n)x^{n-qk},\quad
 c_q(z)=\frac{\prod_{s=0}^{q-1}p(z-sk)}
 {\prod_{t=1}^{q}(\lambda(z)-\lambda(z-tk))},\quad c_0(z)=1,
$$
retaining nonnegative powers of $x$. Here $c_q$ has decrease $qk$
in Lemma~\ref{lem:finite-recurrence-poles}.
Assume a fixed recurrence and $d\ge3$; write
$\lambda(z)=\lambda_dz^d+\lambda_{d-1}z^{d-1}+\cdots$ and
$\omega=e^{2\pi i/d}$. Rescaling $z=qku$ at the simple root
$u=(1-\omega)^{-1}$ of $u^d-(u-1)^d$ gives a root of
$\lambda(z)-\lambda(z-qk)$ with
\refstepcounter{equation}\label{eq:resonance-roots}
$$
 \rho_q=\frac{qk}{1-\omega}
 -\frac{\lambda_{d-1}}{d\lambda_d}+O(q^{-1}).
 \leqno\textup{(\theequation)}
$$
Thus $\operatorname{Im}\rho_q=(qk/2)\cot(\pi/d)+O(1)\to+\infty$.
The numerator zeros have bounded imaginary parts, so $\rho_q$ is
a pole of $c_q$ for large $q$, contradicting
Lemma~\ref{lem:finite-recurrence-poles}. Thus $d\le2$.

Suppose $d=2<\ell$. Complete the missing degrees and take the dual
basis $(\mathbf a_n)_{n\ge0}$. Choose $N=kj_0$ with $N-k$ beyond all
thresholds, spectral collisions and nonnegative integer roots of $p$.
Put $\xi=\lambda(N)$ and $u_i=\duality{\mathbf a_N}{x^{ki}}$. Then
$$
 u_{j_0}=1,\quad
 u_i=-\frac{p(ki)}{\lambda(ki)-\xi}u_{i-1}\ne0,\quad i>j_0.
$$
The transposed recurrence gives
$\boldsymbol\eta=x^k\mathbf a_N
=\sum_{j=N-k}^{N+kw}c_j\mathbf a_j$,
so $\spanC\{(\op M')^t\boldsymbol\eta:t\ge0\}$ is finite-dimensional.
For $\mathbf b\in\Pp'$, its normalized moments
$f(i)=\duality{\mathbf b}{x^{ki}}/u_i$ transform under $\op M'-\xi$, for $i>j_0$, as
$$
 f(i)\longmapsto(\lambda(ki)-\xi)(f(i)-f(i-1)).
$$
For rational $f$ with $\deg_\infty f>0$, the image has degree
$\deg_\infty f+1$.
For $\boldsymbol\eta$, $f(i)=u_{i+1}/u_i$ has degree $\ell-2>0$.
Its successive images have distinct degrees and are linearly independent,
a contradiction.

Conversely, let $d=2$ and $\ell\le2$. Write
$$
 \op{M}=(ax^2+sx+t)\partial^2+(ux+v)\partial+c\op{I},\quad a\ne0.
$$
A direct calculation gives
$$
 \begin{gathered}
 [\op{M},[\op{M},x]]=2a(\op{M}x+x\op{M})+2s\op{M}
 +\beta x+\gamma\op{I},\\[7pt]
 \beta=u(u-2a)-4ac,\quad \gamma=(u-2a)v-2cs.
 \end{gathered}
$$
With $\lambda(n)=an(n+\nu)+c$, $\nu=(u-a)/a$, put
$$
 R_n(z)=(z-\lambda(n))^2-2a(z+\lambda(n))-\beta
 =(z-\lambda(n-1))(z-\lambda(n+1)).
$$
The commutator identity gives
$R_n(\op M)xP_n(x)=(2s\lambda(n)+\gamma)P_n(x)$.
Put $V_n=\Pp_{n+1}$,
$W_n=\spanC\{P_{n-1}(x),P_n(x),P_{n+1}(x)\}$ and
$Q_n(z)=\prod_{j=-1}^{1}(z-\lambda(n+j))$.
Then $Q_n(\op M)xP_n(x)=0$. For large $n$, its three roots occur
once in the triangular diagonal of $\op M|_{V_n}$. In $V_n^\vee$,
$$
 \operatorname{im}\bigl(Q_n(\op M|_{V_n})^\vee\bigr)=W_n^\perp.
$$
Dual separation gives $xP_n(x)\in W_n$, proving the three-term recurrence.

For affine $\lambda$, normalize $\op M$ and apply
Proposition~\ref{prop:commuting-recurrence} for
\eqref{eq:homogeneous-root-criterion}. The commutator expansion of
$e^{-\op A/(ak)}xe^{\op A/(ak)}$ gives
\eqref{eq:homogeneous-recurrence-coefficients};
Proposition~\ref{prop:affine-spectral-algebra} gives minimum order,
even without a fixed recurrence.
\end{proof}

For $\op M=x^j\partial^j+x^a\partial^\ell$, with
$0\le a<\ell$ and $1\le j<\ell$, put
$k=\ell-a$, $\lambda(z)=z^{\underline j}$ and
$p(z)=z^{\underline\ell}$. Theorem~\ref{prop:homogeneous-recurrence}
gives $j=1$ and $k\mid\ell$, recovering the classification in
\cite[Theorem~1.1 and Proposition~5.6]{abst2026}.
When $k\mid\ell$, take
$g(z)=\prod_{r=0}^{\ell/k-1}(z-rk)$ and
$\op G=x^{-1}g(\op E)$; then $\op A=\op G^k$.
For $j=1$, the coefficients, also when $k\nmid\ell$, are given in
\cite[Theorem~5.2]{abst2026}.

\begingroup
\renewcommand{\addcontentsline}[3]{}
\begin{samepage}
\section*{Acknowledgments}

The author thanks B. Shapiro for sharing a copy of \cite{hst2025},
for encouraging him to pursue these questions, and for all he learned
from their joint work \cite{castilloshapiro2026}.
The author is deeply indebted to the late E.~Horozov, whose profound
results on vector orthogonal polynomials and bispectral operators
were fundamental to this work.
\par
\end{samepage}

\begin{samepage}
\section*{Declaration on the Use of Artificial Intelligence}

All mathematical ideas and results presented in this manuscript are the
author's own, developed through traditional mathematical research.
The ideas arose during the preparation of \cite{castilloshapiro2026}
and are applied here to recurrences of fixed length.
Artificial intelligence tools were used to improve the exposition and
replace parts of the author's approach with simpler arguments in
Lemma~\ref{lem:poisson-perturbation} and Lemma~\ref{lem:rational-affine-unit-ladder},
independently check details of the proofs, and search for potential
counterexamples. The author takes full responsibility for the mathematical
content and correctness of the manuscript.
\par
\end{samepage}
\endgroup

\end{document}